  \documentclass[11pt]{amsart}

  \usepackage{etex} 

  \usepackage{blkarray}
  \usepackage{multirow}

  \usepackage{amsmath}
  \usepackage{amssymb}
  \usepackage{amsthm}
  \usepackage{amsfonts}
  \usepackage{booktabs}
  \usepackage[usenames,dvipsnames,svgnames]{xcolor}
  \usepackage{epsfig}
  \usepackage{fancyvrb} 
  \usepackage{hyperref}
  \usepackage{mathrsfs}
  \usepackage{mathtools}
  \usepackage{pictexwd, dcpic}
  \usepackage{enumitem} 

  \usepackage{graphicx}
  
  \usepackage[parfill]{parskip} 
  \usepackage{hyperref,cite}
	\hypersetup{
   	 	colorlinks=true,
   	 	linkcolor=blue,
    		citecolor=red,      
   	 	urlcolor=cyan,
	}
	
\usepackage{imakeidx}
\makeindex

  \usepackage{tikz}
  \usepackage{tikz-cd}
  \usepackage{caption}
  \usetikzlibrary{decorations.markings}

  \newcommand{\FF}{\mathbb{F}}

  \newcommand{\NN}{\mathbb{N}}

  \newcommand{\QQ}{\mathbb{Q}}
  \newcommand{\RR}{\mathbb{R}}

  \newcommand{\ZZ}{\mathbb{Z}}

  \newcommand{\vA}{\mathcal{A}}
  \newcommand{\vB}{\mathcal{B}}
  \newcommand{\vC}{\mathcal{C}}

  \newcommand{\vH}{\mathcal{H}}
  \newcommand{\vI}{\mathcal{I}}

  \newcommand{\vQ}{\mathcal{Q}}
  
  \newcommand{\vS}{\mathcal{S}}

  \newcommand{\Sp}{\operatorname{Sp}}

  \newcommand{\Aut}{\operatorname{Aut}}

  \newcommand{\Comm}{\operatorname{Comm}}

  \newcommand{\Homeo}{\operatorname{Homeo}}
  \renewcommand{\Im}{\operatorname{Im}} 

  \newcommand{\Mod}{\operatorname{Mod}}

  \newcommand{\Ends}{\operatorname{Ends}}

   \newcommand{\PMap}{\operatorname{PMap}}
   \newcommand{\Map}{\operatorname{Map}}

  \newcommand{\co}{\colon\thinspace} 

  \definecolor{lightgrey}{gray}{.85}

  \theoremstyle{definition}
  
  \newtheorem*{rmk}{Remark}
  
  \newtheorem*{ex}{Example}

  \theoremstyle{plain}
  \newtheorem{thm}{Theorem}[section]

  \newtheorem{cor}[thm]{Corollary}
  \newtheorem{prop}[thm]{Proposition}
  \newtheorem{fact}[thm]{Fact}

  \newtheorem{qu}[thm]{Question}
  
  \newtheorem{prob}[thm]{Problem}

  \theoremstyle{definition}

\makeatletter
\newcommand{\mc}{}
\DeclareRobustCommand{\mc}{%
  M%
  \raisebox{\dimexpr\fontcharht\font`M-\height}{%
    \check@mathfonts\fontsize{\sf@size}{0}\selectfont
    c%
  }%
}
  
  \newenvironment{dedication}
{ \vspace{2cm}
   \itshape             
   \raggedleft          
  }

  \title[]{Big mapping class groups: an overview}
  \author[J. Aramayona]{Javier Aramayona}
    \address{ICMAT (CSIC-UAM-UC3M-UCM), Madrid, Spain}
  \email{aramayona@gmail.com}
  \author[N. G. Vlamis]{Nicholas G. Vlamis}
\address{Queens College, City University of New York, Flushing, New York, USA}
\email{nicholas.vlamis@qc.cuny.edu}

\begin{document}
 
   \begin{abstract}
We survey recent developments on mapping class groups of surfaces of infinite topological type. 
  \end{abstract}

  \maketitle


%
     \begin{figure}[h]
 \includegraphics[scale=0.15]{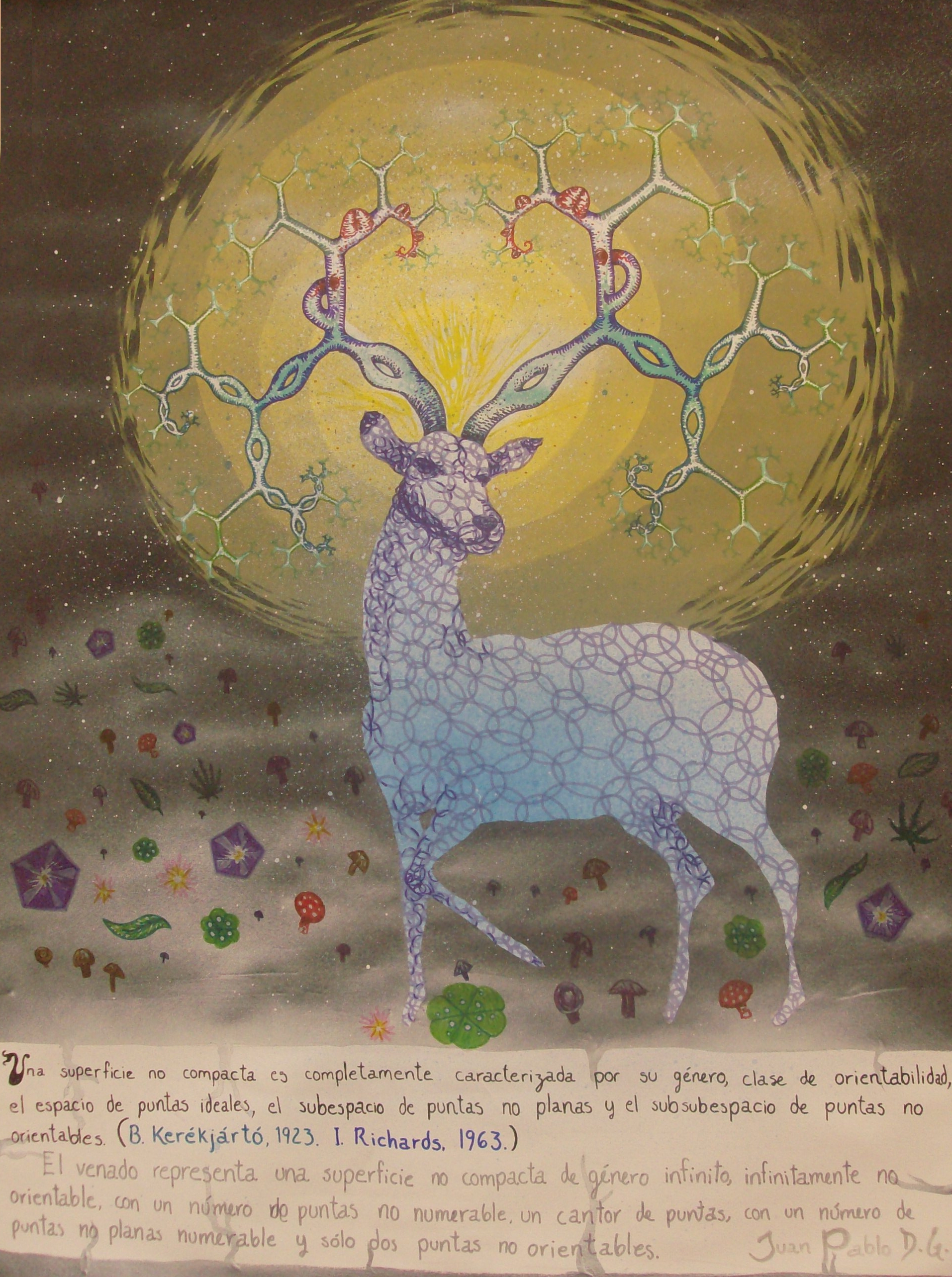} \caption{{\em Infinite-type deer}, by Juan Pablo D\'iaz Gonz\'alez, UNAM}
\end{figure}

      \begin{dedication}
\text{}\\  A Domingo, in memoriam. 
  \end{dedication}

  \section{Introduction}

  In the blogpost  \cite{Cal}, D. Calegari proposed the study of the mapping class group $\Map(\RR^2 \setminus C)$, where $C$ denotes a Cantor set. More concretely, he posed the question of whether this group has an infinite-dimensional space of quasimorphisms, as is the case with the mapping class group of a surface of finite topological type, after a celebrated result of Bestvina--Fujiwara \cite{BeF}. In addition,  Calegari suggested a line of attack on the problem, in analogy with Bestvina--Fujiwara's original argument; in a nutshell, the first idea is to prove that a certain {\em complex of arcs} on which $\Map(\RR^2 \setminus C)$ acts is hyperbolic and has infinite diameter, and then exhibit elements which act {\em weakly properly discontinously} \cite{BeF} on this complex. 
  
This strategy was successfully implemented by J. Bavard in her thesis \cite{Bav} (English translation: \cite{Bav2}), and has since caused a surge of interest in mapping class groups of infinite-type surfaces (or {\em big mapping class groups}\index{mapping class group!big}, in the terminology coined by Calegari) among the geometric group theory and low-dimensional topology communities. Most of the results to date have focused on the basic structure of big mapping class groups, as well as on the similarities and differences with mapping class groups of finite-type surfaces.

This said, big mapping class groups made their appearance in other related areas of mathematics quite a long time ago. For instance, 
big mapping class groups arise naturally in the context of stable properties of mapping class groups \cite{Miller}; infinite-type surfaces are intimately related to the study of quasiconformal maps \cite{Bers}; the so-called {\em braided Thompson's group} $BV$ of Brin \cite{Brin} and Dehornoy \cite{Dehornoy} is naturally a subgroup of the mapping class group of a sphere minus a Cantor set; etc.

The aim of this survey is to give an overview of the recent developments around big mapping class groups, mainly from the point of view of geometric group theory, and to describe some of the connections to other areas of mathematics, such as Polish groups and Thompson's groups. Along the way, we will offer open problems related to the topics covered. 

\medskip

\noindent{\bf Plan of the chapter.} All the objects and definitions needed in the exposition are introduced in Section \ref{sec:prelim}. In Section \ref{sec:important}, we present two results which are crucial to a large number of the results discussed in subsequent sections. Section \ref{sec:topology} deals with topological aspects of big mapping class groups: generation, Polish structure, etc. Section \ref{sec:algebra} concerns algebraic results: automorphisms, homology, relation with Thompson's groups, etc. Finally, in Section \ref{sec:geometry} we will concentrate on the action of big mapping class groups on various hyperbolic complexes constructed from arcs and/or curves on the surface. 

\medskip

\noindent{\bf Big absences.} There are a number of interesting topics related to big mapping class groups which are not covered in this survey. Notably, the relation between mapping class groups and  dynamics \cite{Calcxdyn}, the  theory of Teichm\"uller spaces of infinite-type surfaces (see \cite{Matsuzaki,LiPa} and the references therein), and the theory of infinite translation surfaces (see for instance \cite{Randecker} and the references therein). 

\bigskip

\noindent{\bf Acknowledgements.} We would like to thank Athanase Papadopoulos for inviting us to write this survey. 

 This article is heavily influenced by the AIM workshop ``Surfaces of infinite type'' (29 April--3 May 2019). We would like to thank the organizers of the workshop, J. Bavard, A. Randecker, P. Patel, and J. Tao; we are also grateful to AIM for its hospitality and financial support.
 
J.A. was supported by grants RYC-2013-13008 and PGC2018-101179-B-I00. He acknowledges financial support from the Spanish Ministry of Science and Innovation, through the ``Severo Ochoa Programme for Centres of Excellence in R\&D'' (SEV-2015-0554) and from the Spanish National Research Council, through the ``Ayuda extraordinaria a Centros de Excelencia Severo Ochoa'' (20205CEX001).

N.G.V. was partially supported by PSC-CUNY grant 62571-00 50.
 
Finally, J.A. is indebted to his baby daughter Berta for improving his ability at typing with one hand.

\section{Preliminaries}
\label{sec:prelim}

In this section we introduce the background material needed for the rest of the article.

\subsection{Surfaces and their classification} Throughout this article, all surfaces\index{surface} considered will be assumed to be second countable, connected, orientable, and have compact (possibly empty) boundary. If the fundamental group of $S$ is finitely generated, we will say that $S$ is of {\em finite type} \index{topological type}\index{topological type!finite}; otherwise, we will say that $S$ is of {\em infinite type}\index{topological type!infinite}.

It is well-known that the homeomorphism type of a finite-type surface is determined by the triple $(g,p,b)$, where $g \ge 0$ is the genus, and $p,b\ge 0$ are, respectively, the number of punctures and boundary components of the surface. Because of this fact, we will use the standard notation $S_{g,p}^b$ for  the surface specified by these data; as usual, we will drop $p$ and $b$ from the notation whenever they are equal to zero. 

There is also a similar classification for infinite-type surfaces \cite{Ke, Richards}, in terms of genus, number of boundary components, and the topology of the {\em space of ends}, which we now define. First, an {\em exiting sequence} is a sequence\index{exiting sequence} $\{U_n\}_{n\in \NN}$ of connected open subsets of $S$ with the following properties: 

\begin{enumerate}
\item $U_n \subset U_m$ whenever $m<n$, 
\item $U_n$ is not relatively compact for any $n \in \NN$, 
\item $U_n$ has compact boundary for all $n \in \NN$, and
\item any relatively compact subset of $S$ is disjoint from all but finitely many $U_n$'s. 
\end{enumerate}

Two exiting sequences are {\em equivalent} if every element of the first is eventually contained in some element of the second, and vice versa. We denote by $\Ends(S)$ the set of all equivalence classes of exiting sequences of $S$; an element of \( \Ends(S) \) is referred to as an \emph{end}\index{end}\index{end!topological} of \( S \). 
The set \( \Ends(S) \) becomes a topological space, called \emph{the space of ends of \( S \)}\index{space of ends}\index{end!space of}, by specifying the following basis: given a subset $U \subset S$ with compact boundary, consider the set $U^*$ of all ends represented by an exiting sequence eventually contained in $U$; the set \( \{ U^* : U \subset S \text{ open with compact boundary}\} \) is the desired basis.  
If \( U \) is an open set with compact boundary and \( e \in  U^* \), then we say that \( U \) is a \emph{neighborhood} of the end \( e \). 

Given the above basis, it is not difficult to see that \( \Ends(S) \) is Hausdorff, totally disconnected, and second countable. 
Moreover, the definition above can be reframed to describe \( \Ends(S) \) in terms of an inverse limit of compact spaces; in particular, Tychonoff's theorem implies \( \Ends(S) \) is compact.
(For a reference, see \cite[Chapter 1]{ASario}.) 

\begin{thm}
For any surface $S$, the space $\Ends(S)$ is totally disconnected, second countable, and compact. In particular, \( \Ends(S) \) is homeomorphic to a closed subset of a Cantor set.
\label{thm:endsdisc}
\end{thm}

We now proceed to describe the classification of infinite-type surfaces up to homeomorphism. 
To this end, we will say that an end is \emph{planar}\index{end!planar} if it admits a neighborhood that is embeddable in the plane; otherwise an end is \emph{non-planar}\index{end!non-planar} (or \emph{accumulated by genus}\index{end!accumulated by genus}) and every neighborhood of the end has infinite genus.
Denote by $\Ends_{np}(S)$ the subspace of $\Ends(S)$ consisting of non-planar ends, noting that it is closed in the subspace topology. The following result was proved by Ker\'ekj\'art\'o \cite{Ke} and Richards \cite{Richards}. 

\begin{thm}[Classification\index{surface!classification of}\index{classification of surfaces}, \cite{Ke,Richards}]
Let $S_1,S_2$ be surfaces, and write $g_i$ and $b_i$, respectively, for the genus and number of boundary components of $S_i$. Then $S_1 \cong S_2$ if and only if $g_1= g_2$, $b_1= b_2$ and there is a homeomorphism \[\Ends(S_1) \to  \Ends(S_2)\] that restricts to a homeomorphism \[\Ends_{np}(S_1) \to \Ends_{np}(S_2).\]
\end{thm}

In light of the above result, an obvious question is: given two closed subsets $X,Y$ of a Cantor set, with $Y\subset X$, can they be realized as the spaces of ends (resp. ends accumulated by genus) of some surface? The following theorem, due to Richards \cite{Richards}, states that the answer is ``yes": 

\begin{thm}[Realization, \cite{Richards}]
Let $X, Y$ be closed subsets of a Cantor set with \( Y \subset X \). Then there exists a surface $S$ such that $\Ends(S)\cong X$ and $\Ends_{np}(S) \cong Y$. 
\end{thm} 

With the classification and realization theorems at hand, we make a quick note about cardinality:
there are exactly $\aleph_0$ many homeomorphism classes of compact surfaces, but $2^{\aleph_0}$ many homeomorphism classes of second-countable surfaces. 
The second statement follows from a count on the homeomorphism classes of closed subsets of the Cantor set \cite{Reichbach}.
Interestingly, if one drops the condition of second countability, then there are $2^{\aleph_1}$ many homeomorphism classes of surfaces \cite{Gauld}.

\subsubsection{Some important examples} Several infinite-type surfaces have standard names, which makes them easy to identify; these are as follows: 

\begin{figure}[t]
\includegraphics[width=0.8\textwidth]{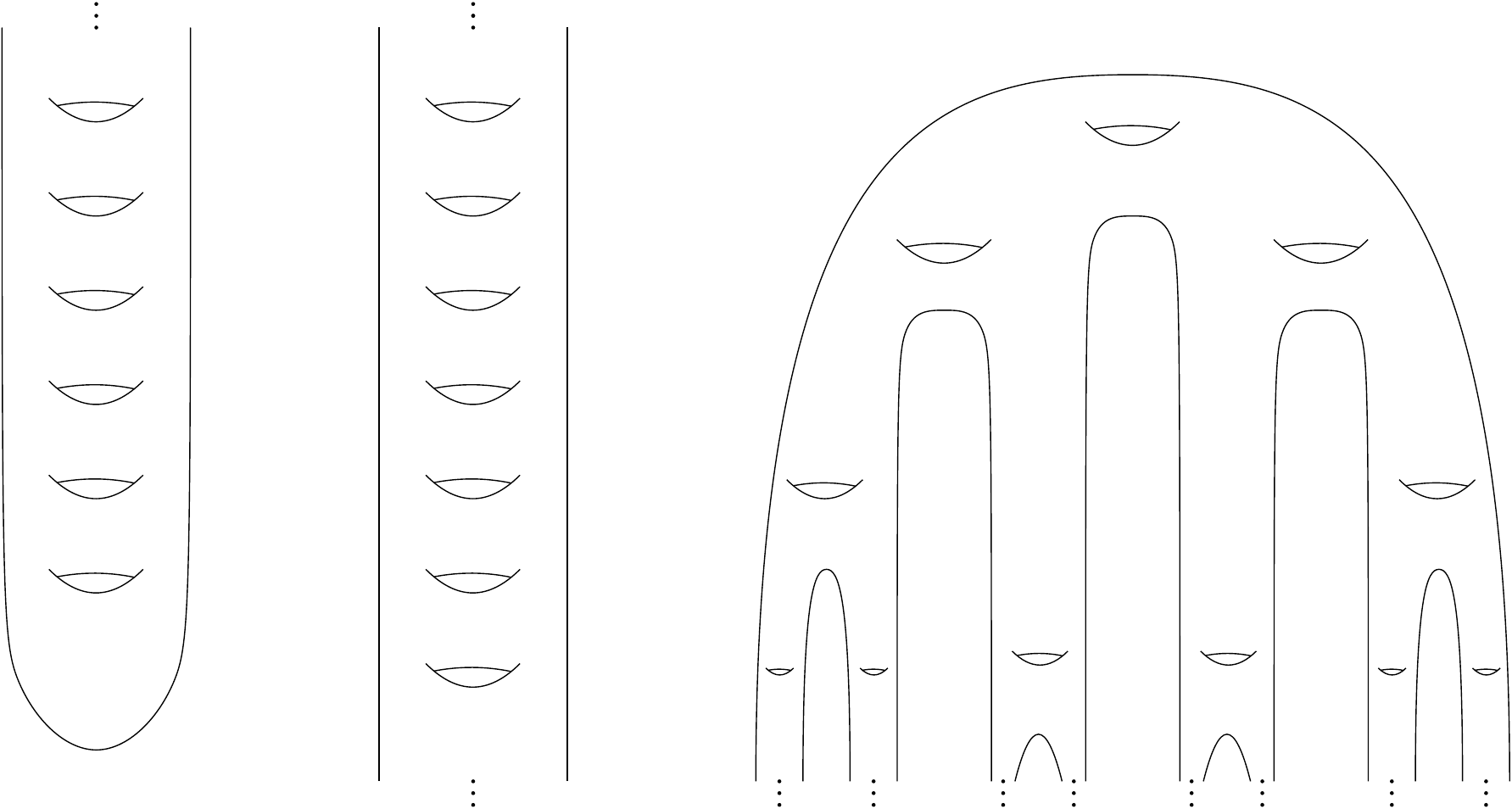}
\caption{From left to right: Loch Ness monster surface, Jacob's ladder surface, and the blooming Cantor tree surface.}
\label{fig:surfaces}
\end{figure}

\begin{itemize}
\item {\em The Loch Ness monster surface\index{sufrace!Loch Ness monster}}: the infinite-genus surface with exactly one end (which is necessarily non-planar). 
\item {\em Jacob's ladder surface\index{surface!Jacob's ladder}}: the infinite-genus surface with exactly two ends, both non-planar. 
\item {\em The Cantor tree surface\index{surface!Cantor tree}}: the planar surface whose space of ends is a Cantor space. Hence, this surface is homeomorphic to a sphere minus a Cantor set. 
\item {\em The blooming Cantor tree surface\index{surface!blooming Cantor tree}}: the infinite-genus surface whose space of ends is a Cantor space, and such that every end is non-planar.
\item {\em The flute surface\index{surface!flute}}: the planar surface whose space of ends has a unique accumulation point.  Hence, this surface is homeomorphic to \( \mathbb C \smallsetminus \mathbb Z \) (and the end space is homeomorphic to \( \{0\} \cup \{ \frac1n : n\in\NN\} \), viewed as a subset of \( \RR \)).
\end{itemize}

The Loch Ness monster surface, Jacob's ladder surface, and the blooming Cantor tree surface are shown in Figure \ref{fig:surfaces}; the Cantor tree surface can be seen in Figure \ref{fig:genus0}.
To the authors' knowledge, the first two of these names were introduced by Phillips--Sullivan \cite{Phillips}, the second two by Ghys \cite{Ghys2}, and the last by Basmajian \cite{BasmajianCollar}.
It is worth noting that in \cite{Ghys2}, Ghys shows that a generic non-compact leaf of 2-dimensional lamination of a metric space is either the plane, the cylinder, or one of the first four surfaces above.

\subsection{Arcs and curves}
By an {\em arc}\index{arc} on $S$ we mean the homotopy class of a properly embedded copy of $\RR$. 
Abusing notation, we will not distinguish between arcs and their representatives. Two arcs are {\em disjoint} if they have disjoint representatives; otherwise we say that they intersect. The {\em intersection number}\index{intersection number}, denoted \( i(\cdot, \cdot) \), between two arcs is the minimum (possibly infinite) number of points of intersection between representatives.

By a {\em curve}\index{curve} on $S$ we mean the homotopy class of a simple closed curve on $S$ which does not bound a disk, a punctured disk, or an annulus whose other boundary component is contained in $\partial S$. As was the case with arcs, we will use the same notation for curves and their representatives. We say that a curve $\alpha$ is {\em non-separating}\index{curve!non-separating} if $S \smallsetminus \alpha$ is connected; otherwise we say that $\alpha$ is {\em separating}\index{curve!separating}. 
Again, we may talk about when two curves are disjoint or intersect, and define their intersection number as we did with arcs and use the same notation. Note, however, that the intersection number between two curves is necessarily a finite number. 

A {\em multicurve}\index{curve!multi-} is a set of pairwise-distinct and pairwise-disjoint curves. A {\em pants decomposition}\index{pants!decomposition} is a multicurve $P$ that is maximal with respect to inclusion, and such that any compact set on $S$ is intersected by only finitely many elements of $P$. As such, the interior of every connected component of the complement of $P$ in $S$ is homeomorphic to a sphere with three points removed, commonly referred to as a {\em pair of pants}\index{pants!pair of}.

\subsection{Mapping class group}\index{mapping class group} Consider the group $\Homeo(S, \partial S)$ of homeomorphisms of $S$ that restrict to the identity on the boundary of $S$, equipped with the compact-open topology\index{mapping class group!topology of}\index{topology!compact-open}, and the subgroup $\Homeo^+(S,\partial S)$ consisting of those elements that preserve orientation.  Let $\Homeo_0(S,\partial S)$ denote the path component of the identity in $\Homeo(S,\partial S)$, and note that \( \Homeo_0(S, \partial S) \subset \Homeo^+(S, \partial S) \). 
The {\em extended mapping class group}\index{mapping class group!extended} is 
\[
\Map^\pm(S):= \Homeo(S,\partial S) / \Homeo_0(S,\partial S),
\] 
and the  {\em mapping class group} is the subgroup 
\[
\Map(S):= \Homeo^+(S,\partial S) / \Homeo_0(S,\partial S).
\] 
The extended mapping class group becomes a topological group with the quotient topology coming from the compact-open topology on $\Homeo(S,\partial S)$.
Combining \cite[Theorem 6.4]{Epstein} and \cite[Theorem 1]{Fox}, we see that the elements of \( \Map(S) \) are exactly the isotopy classes of orientation-preserving homeomorphisms of \( S \) (see the appendix in \cite{vlamisnotes} for a more detailed discussion). 

(Note that is not clear or obvious that the mapping class group is Hausdorff, since---a priori---path components are not closed subsets. Being Hausdorff is a condition that is often required in the definition of topological group.  We will deal with this in Section \ref{sec:topology}.)

\subsection{Several natural subgroups}
Throughout the survey, several natural subgroups of mapping class groups will appear: we provide their definition here. 

\subsubsection{Pure mapping class group}\index{mapping class group!pure}
Observe that every homeomorphism of $S$ induces a type-preserving homeomorphism of its space of ends.  In
other words, there is a natural map 
\begin{equation}
\Homeo^+(S,\partial S) \to \Homeo(\Ends(S), \Ends_{np}(S)),
\label{eq:surjhom}
\end{equation} where the latter group denotes the subgroup of the homeomorphism group of $\Ends(S)$ whose elements preserve $\Ends_{np}(S)$ setwise. One checks this is a continuous homomorphism when  $\Homeo(\Ends(S), \Ends_{np}(S))$  is equipped with the (subgroup topology coming from the) compact-open topology. 

Richards's proof of the classification of surfaces can readily be adapted to establish the surjectivity of the homomorphism given in \eqref{eq:surjhom}.
As an isotopy fixes every end of a surface, the homomorphism \eqref{eq:surjhom} factors through \( \Map(S) \) yielding a surjective homomorphism
\begin{equation}\Map(S) \to \Homeo(\Ends(S), \Ends_{np}(S)).
\end{equation}
The {\em pure mapping class group}, written \( \PMap(S) \), is the kernel of the above homomorphism. In particular, we have a short exact sequence \begin{equation}
1 \to \PMap(S) \to \Map(S) \to \Homeo(\Ends(S), \Ends_{np}(S)) \to 1
\label{eq:SEShomeo}
\end{equation}

It is worth noting that by Stone's representation theorem, there is a one-to-one correspondence (or, technically, a contravariant functor) between closed subsets of the Cantor set and countable Boolean algebras.
There is a large amount of literature about automorphism groups of boolean algebras, which can be translated to homeomorphism groups of end spaces of surfaces (and vice versa).

We also note that, by the  definition of the mapping class groups, \( \Map(S) = \PMap(S) \) if and only if either  \( |\Ends(S)| \leq 1 \) or \( |\Ends(S)| = 2 \) and \( S \) has exactly  one planar end.

\subsubsection{Compactly supported mapping class group}\index{mapping class group!compactly supported} 
An element of \( \Map(S) \) is \emph{compactly supported} if it has a representative homeomorphism that is the identity outside of a compact subset. 
The \emph{compactly supported mapping class group}, denoted \( \Map_c(S) \),   
is the subgroup of $\Map(S)$ consisting of the compactly supported elements.
Observe that, in fact, $\Map_c(S) < \PMap(S)$. 

We say a compact subsurface \( X \) of a surface \( S \) is \emph{essential} if no component of \( S\smallsetminus X \) is a disk or annulus.
If $X$ is an essential compact subsurface of $S$, then $\Map(X) < \Map_c(S)$.
Note that for any two essential compact subsurfaces \( X \) and \( Y \) of \( S \), we have  \( \Map(X) < \Map(Y) \) whenever \( X \subset Y \). 
Moreover, the union of all compact subsurfaces of \( S \) is equal to \( S \); hence, we have:

\begin{prop}
\label{prop:gencompact}
For any surface $S$, \[\Map_c(S) = \lim_{\to} \Map(X),\] where the direct limit is taken over all essential compact subsurfaces \( X \) of $S$, ordered by inclusion.
\end{prop}

\subsubsection{Torelli group}\index{Torelli group}
Observe that every element of $\Map(S)$ acts on the homology group $H_1(S,\ZZ)$ by automorphisms.  In other words, there is a homomorphism {\begin{equation}
\Map(S) \to \Aut(H_1(S,\ZZ)).
\label{eq:homrep}
\end{equation}
We remark that if $S$ is a finite-type surface of genus $g$ and with at most one puncture, then $\Aut(H_1(S,\ZZ))$ is isomorphic to the symplectic group $\Sp(2g,\ZZ)$, although this is not true in general. The {\em Torelli group} $\vI(S)$ is the kernel of the homomorphism \eqref{eq:homrep}; in other words, it is the subgroup of $\Map(S)$ whose elements act trivially on homology. 
Observe that \( \vI(S) \) is a subgroup of \( \PMap(S) \).

\subsection{Modular groups}\index{modular group}

Naturally associated to a Riemann surface is the subgroup \( \rm{QC}(X) \) of \( \Homeo^+(X) \) consisting of the quasi-conformal homeomorphisms.
The image of \( \rm{QC}(X) \) in \( \Map(X) \), denoted \( \Mod(X) \),  is commonly referred to as either the \emph{Teichm\"uller modular group} of \( X \) or the \emph{quasi-conformal mapping class group}\index{mapping class group!quasi-conformal} of \( X \). 
In the case that \( X \) is of finite topological type, \( \Mod(X) \) and \( \Map(X) \) agree and are routinely interchanged in the literature; however, this fails to be the case for infinite-type surfaces.

In the infinite-type setting, unlike mapping class groups, modular groups have a long history of being studied, especially from the theory of Riemann surfaces and Teichm\"uller theory. 
As such, discussing the modular group would be a survey in-of-itself and we will make no further mention of it.
But, we note that there are surely many interesting questions and problems related to how \( \Mod(X) \) sits as a subgroup of \( \Map(S) \), where \( X \) is a Riemann surface homeomorphic to an infinite-type surface \( S \). 


\section{Two important results}
\label{sec:important}

In this section we present two results that underpin a large number of the topics discussed in latter sections. Throughout this section, every surface is assumed to have empty boundary.

\subsection{Alexander method}\index{Alexander method}

As mentioned in the introduction, \( \Map(S) \) inherits a natural topology when viewed as a quotient of \( \Homeo^+(S) \), equipped with the compact-open topology.
It is standard to require that a topological group be Hausdorff, and so it is not immediately obvious that \( \Map(S) \) in this topology is in fact a topological group. 
However, we can use the extension of Alexander's method to infinite-type surfaces given in \cite{HMV}.
Here, we state the corollary we require:

\begin{thm}[{\cite[Corollary 1.2]{HMV}}]
\label{thm:hausdorff}
Let \( S \) be an infinite-type surface. 
If \( f \in \Homeo^+(S) \) fixes the isotopy class of every simple closed curve, then \( f \) is isotopic to the identity. 
\end{thm}

Theorem \ref{thm:hausdorff} can used  to separate the identity from any other element in \( \Map(S) \) by an open set and, for topological groups, this is enough to guarantee the group is Hausdorff; hence, \( \Map(S) \) is a topological group.

\subsection{Automorphisms of the curve graph}\index{curve graph!automorphisms of} The {\em curve graph}\index{curve graph} $\vC(S)$ of $S$ is the simplicial graph whose vertex set is the set of curves on $S$, and where two vertices are adjacent in $\vC(S)$ if and only if the corresponding curves on $S$ are disjoint. From now on we will not distinguish between vertices of $\vC(S)$ and the curves they represent. 

Observe that $\Map^\pm(S)$ acts on $\vC(S)$ by simplicial automorphisms. In fact, the combined work of Ivanov \cite{Iva}, Korkmaz \cite{Kor}, and Luo \cite{Luo} shows that, with the exception of the twice-holed torus, there are no other automorphisms of \( \vC(S) \) when \( S \) is of finite type. In the infinite-type setting, the analogous result was proved independently by Hern\'andez--Morales--Valdez \cite{HMV-iso} and Bavard--Dowdall--Rafi \cite{BDR}: 

\begin{thm}
\label{thm:ivanovcurves}
If \( S \) is an infinite-type surface, then the group of simplicial automorphisms of $\vC(S)$ is naturally isomorphic to $\Map^\pm(S)$. 
\end{thm} 

Note that, in particular, Theorem \ref{thm:hausdorff} is required to show that the action of \( \Map(S) \) on \( \vC(S) \) has no kernel.



\section{Topological aspects} 
\label{sec:topology}

We will see in this section that big mapping class groups are interesting topological groups---a divergence from the finite-type setting.
This offers exciting new connections for mapping class groups, some of which we explore below.

It follows from the Alexander method for finite-type surfaces (see \cite[Proposition 2.8]{FM}) that \( \Map(S) \) is discrete when \( S \) is of finite-type.
However, this is far from true for big mapping class groups: to see this, let \( S \) be an infinite-type surface and let \( \{c_n\}_{n\in\NN} \) be a sequence of simple closed curves such that, for every compact subset \( K \) of \( S \), there is an integer \( N \) such that \( K \cap c_n = \emptyset \) for all \( n > N \). 
If \( T_n \) is the Dehn twist about \( c_n \), then the sequence \( \{ T_n \}_{n\in\NN} \) limits to the identity in \( \Map(S) \). 

\subsection{The permutation topology}\index{topology!permutation}\index{mapping class group!topology of}

In order to investigate the topology of \( \Map(S) \) in more depth, it is convenient to have a more combinatorial description of its topology. 

Let $\Gamma$ be a simplicial graph with a countable set of vertices, and let $\Aut(\Gamma)$ be the group of simplicial automorphisms of  $\Gamma$. 
Given a subset $A$ of $\Gamma$, let 
\[U(A):= \{g \in \Aut(\Gamma) \mid g(a) = a \text{ for all } a \in A\}.\]
Then $\Aut(\Gamma)$ may be endowed with a natural topology, called the {\em permutation topology}, defined by declaring the $\Aut(\Gamma)$-translates of $U(A)$, for every finite subset $A$ of $\Gamma$, a basis for the topology.
Equivalently, the permutation topology is the coarsest topology in which, for every \( v \in \vC(S) \), the function \( \omega_v \co \Aut(\Gamma) \to \Gamma \) defined by \( \omega_v(g) = g(v) \) is continuous.

With respect to the permutation topology, $\Aut(\Gamma)$ becomes a second countable (and in particular, separable) topological group. 
Moreover, it is a standard exercise in descriptive set theory texts to show that \( \Aut(\Gamma) \) supports a complete metric (which---usually---fails to be \( \Aut(\Gamma) \)-invariant).

In particular, \( \Aut(\Gamma) \) is an example of a \emph{Polish group}\index{Polish group}\index{group!Polish}, that is, a separable and completely metrizable group.
Polish groups are a well-studied class of groups and we will make use of their theory.

For an infinite-type surface \( S \) with empty boundary, let \( \Gamma = \mathcal  C(S) \), then, by Theorem \ref{thm:ivanovcurves}, we can identify \( \Map^\pm(S) \) with \( \Aut(\Gamma) \) and equip \( \Map^\pm(S) \) with the associated permutation topology.
It is an exercise in definitions and the Alexander method to show that this permutation topology agrees with the compact-open topology. 
Recall that a \( G_\delta \) subset of a topological space is a subset that can be written as the intersection of countably many open sets (note that in a metrizable space, every closed set is a \( G_\delta \) subset). 
As a consequence of this discussion, we have:

\begin{prop}
\label{prop:polish}
Let $S$ be a infinite-type surface, possibly with non-empty boundary. Then, \( \Map^\pm(S) \) and all its \( G_\delta \)-subgroups, including \( \Map(S) \) and \( \PMap(S) \), are Polish. 
\end{prop} 

Note that, unlike the preceding discussion, Proposition \ref{prop:polish} does not require \( S \) to have empty boundary: this is because the mapping class group of a bordered surface can be embedded in a borderless surface as a closed subgroup.

\subsection{Basic properties}

Now that we have an understanding of the topology of mapping class groups, we can investigate their basic properties.  
First, note that the sets in the basis defined above for \( \Map(S) \) are in fact clopen and hence mapping class groups are zero-dimensional.  

Now, let \( S \) be of infinite type. 
Observe that if \( A \subset \mathcal C(S) \) and \( c \in \vC(S) \) such that \( c \cap a = \varnothing \) for all \( a \in A \), then the sequence \( \{T_c^n\}_{n\in\NN} \) has no limit point and is contained in \( U(A) \); in particular, again by homogeneity, we can conclude that every compact subset of \( \Map(S) \) is nowhere dense.
This also establishes the weaker fact that \( \Map(S) \) fails to be locally compact.
Moreover, as a Polish space cannot be the countable union of nowhere dense subsets, we can conclude that \( \Map(S) \) is not compactly generated.\footnote{There are two standard meanings for compactly generated, one algebraic and one toplogical.  For clarity, we are referring to the algebraic setting: specifically, we mean that if a set \( \mathcal S \) generates \( \Map(S) \), as a group, then \( \mathcal S \) cannot be compact.}
Lastly, the Alexandrov--Urysohn Theorem (see \cite[Theorem 7.7]{Kechris}) establishes \( \NN^\NN \) as the unique space, up to homeomorphism, that is non-empty, Polish, zero-dimensional, and in which every compact subset has non-empty interior; hence, \( \Map(S) \) is homeomorphic to \( \NN^\NN \). 
We record these observations in the following theorem:

\begin{thm}\index{mapping class group!topology of}
\label{thm:basics}
For every infinite-type surface \( S \), 
\begin{enumerate}
\item \( \Map(S) \) is not locally compact,
\item \( \Map(S) \) is not compactly generated,
\item \( \Map(S) \) is homeomorphic to the Baire space \( \NN^\NN \) (which in turn is homeomorphic to \( \RR \smallsetminus \QQ) \). 
\end{enumerate}
\end{thm}

Theorem \ref{thm:basics} establishes big mapping class groups as large topological groups. 
It is often the case that the tools developed for studying finitely-generated groups have natural analogs in the setting of locally-compact compactly-generated topological groups.
The failure of big mapping class groups to fall into this category will generally complicate matters, but simultaneously offers big mapping class groups as potential fertile ground for applying the rapidly developing and exciting theory and tools of non-locally-compact topological groups.
We will see this below when we discuss the geometry of mapping class groups. 

\subsection{Topological generation}

Since big mapping class groups are separable, they are necessarily topologically generated by a countable set, that is, there exists a countable set that generates a dense subgroup.
The goal of this subsection is to produce such a topological generating set whose elements are relatively simple. 
Recall that for a connected finite-type surface \( S \), its pure mapping class group \( \PMap(S) \) is generated by---a finite set of---Dehn twists.
In order to generate the full mapping class group, it is necessary to add half-twists, which correspond to transpositions in the symmetric group isomorphic to \( \Map(S) / \PMap(S) \).  

In the infinite-type setting, equation \eqref{eq:SEShomeo} tells us that \( \Map(S) / \PMap(S) \) is isomorphic to \( \Homeo(\Ends(S), \Ends_{np}(S)) \), so in order to understand topological generating sets for \( \Map(S) \), we would also have to do so for the latter homeomorphism groups; this will take us too far afield and so we will focus on generating \( \PMap(S) \). 

Using the fact that the mapping class group of a compact surface is generated by Dehn twists, we see that the group \( \Map_c(S) \) consisting of compactly supported mapping classes is generated by Dehn twists. 
It is natural to ask if the closure of this group is all of \( \PMap(S) \). The next result, proved by Patel and the second author in \cite{PV}, shows that this is true only in certain cases: 

\begin{thm}(\cite{PV})
\label{thm:PV}
The set of Dehn twists topologically generate \( \PMap(S) \) if and only if \( S \) has at most one non-planar end.
\index{mapping class group!topologically generating}
\end{thm}

The only impediment to Dehn twists topologically generating is the existence of a homeomorphism \( f \co S \to S \) and a separating curve \( \gamma \) non-trivial in homology such that \( f(\gamma) \cap \gamma = \varnothing \).  
As it turns out, this can only be done---while fixing the ends---if there are at least two non-planar ends.
Let us give an example of such a homeomorphism, known as a \emph{handle shift}, which was introduced in \cite{PV}.

For \( n \in \ZZ \), let \( B^\pm_n \) be the open Euclidean disks of radius 1 in \( \RR^2 \) centered at \( (n, \pm 2) \), respectively.
Let \( \Sigma \) be the (infinite-genus) surface obtained from \( \RR \times [-4,4] \) by, for each \( n \in \ZZ \), removing \( B^\pm_n \) and identifying \( \partial B^+_n \) and \( \partial B^-_n \) via an orientation-reversing homeomorphism. 
Up to isotopy, there is a unique homeomorphism \( \sigma \co \Sigma \to \Sigma \) determined by requiring

\begin{enumerate}
\item \( \sigma((x,y)) = (x+1, y) \) for all \( (x,y) \in \Sigma \) with \( |y| \leq 3  \), and
\item \( \sigma((x,\pm4)) = (x, \pm 4) \) for all \( x \in \RR \). 
\end{enumerate}

\begin{figure}[t]
\includegraphics{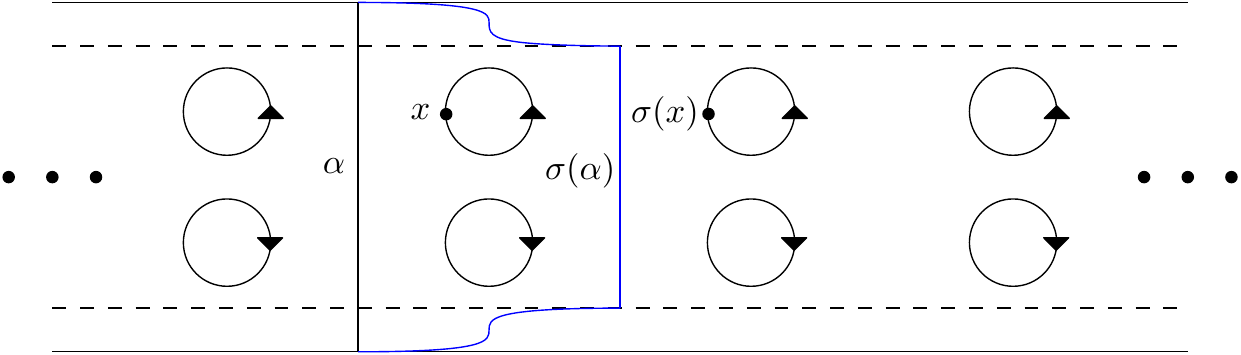}
\caption{The circles are identified vertically to obtain \( \Sigma \).}
\label{fig:handleshift}
\end{figure}

See Figure \ref{fig:handleshift} to see the behavior of \( \sigma \) on a vertical arc. 
Now, for an infinite-genus surface \( S \), we say a homeomorphism \( h \co S \to S \) is a \emph{handle shift}\index{handle shift} if there exists a proper embedding \( \iota \co \Sigma \to S \) such that 

\[
h = \left\{ 
\begin{array}{ll}
\iota \circ \sigma \circ \iota^{-1}(x) & x \in \iota(\Sigma) \\
x & \text{otherwise}
\end{array} \right.
\]

We will also refer to a mapping class containing a handle shift as a handle shift itself. 
Identifying \( \Sigma \) with its image under \( \iota \), we say that \( h \) is \emph{supported} on \( \Sigma \).
Since the embedding \( \iota \) is required to be proper, there is an induced map \( \iota_\infty \co \Ends(\Sigma) \to \Ends(S) \).
It follows that \( h \) has an \emph{attracting} and a \emph{repelling} end, which we label \( h^+ \) and \( h^- \) respectively, and that satisfy \[ \lim_{n \to \pm \infty} h^n(x) = h^\pm \]
for every \( x \) in the interior of \( \Sigma \) (the limit is formally taken in the Freudenthal compactification of \( S \)).
Note that if \( h_1 \) and \( h_2 \) are isotopic handle shifts, then \( h_1^\pm = h_2^\pm \); therefore, we can talk about the attracting and repelling ends of a mapping class associated to a handle shift.  

Let \( h \) be a handle shift supported on \( \Sigma \) in an infinite-genus surface \( S \) with at least two non-planar ends and such that \( h^+ \neq h^- \). 
Now observe that if we take a separating curve \( \gamma \) that is non-trivial in homology and such that \( \gamma \cap \Sigma \) is connected and isotopic to a vertical arc, then \( \gamma \) is non-trivial in homology, \( \gamma \) is not homotopic to \( h(\gamma) \), and \( i(\gamma, h(\gamma)) = 0 \).
As described in \cite{PV}, these conditions guarantee that \( h \) is not a limit of compactly supported mapping classes.

It was shown in \cite{PV} that the set of Dehn twists together with the set of handle shifts topologically generate \( \PMap(S) \). 
But, the set of handle shifts is uncountable and we want a countable dense subset.
As a corollary of a---much stronger---result in \cite{APV}, we can reduce to a countable collection:

\begin{thm}(\cite{APV})
\label{thm:topgen}
If \( S \) is an infinite-genus surface with at least two non-planar ends, then there exists a countable set consisting of Dehn twists and handle shifts topologically generating \( \PMap(S) \). 
\end{thm}

The handle shifts obtained from \cite{APV} will pairwise commute; however, for a weaker, but direct version, it would suffice to choose a countable dense subset \( \{(e^+_n, e^-_n)\}_{n\in\NN} \) in \( \Ends_{np}(S) \times \Ends_{np}(S) \) and handle shifts \( h_n \in \PMap(S) \) such that \( h_n^\pm = e^\pm_n \). 
It can be checked that these handle shifts along with Dehn twists will topologically generate \( \PMap(S) \). 

Adapting an argument presented in \cite[Theorem 7.16]{FM} showing that the mapping class group of a finite-type surface is generated by torsion elements, Afton--Freedman--Lanier--Yin \cite{AFLY} observed: 

\begin{thm}[\cite{AFLY}]
If $S$ is an infinite-genus surface, then $\PMap(S)$ is topologically generated by handle shifts. 
\end{thm}

\subsubsection{Torelli group}\index{Torelli group!topologically generating}

As noted previously, \( \vI(S) \) is contained in \( \PMap(S) \); moreover, handle shifts act non-trivially on homology and hence \( \vI(S) \) contains no handle shifts.
This is enough to imply that  \( \vI(S) \) is contained in the closure of \( \Map_c(S) \) (this follows from Corollary \ref{cor:semidirect} below).
Letting \( \vI_c(S) \) denote the intersection \( \vI(S) \cap \Map_c(S) \), it is natural to ask if the closure of \( \vI_c(S) \) is all of \( \vI(S) \). The answer is yes:

\begin{thm}[\cite{ATorelli}]
If \( S \) is an infinite-type surface, then \( \vI_c(S) \) is dense in \( \vI(S) \). 
\end{thm}

Combining results of Birman \cite{Birman}, Powell \cite{Powell} and an argument due to Justin Malestein, the above theorem implies the following (see \cite{ATorelli} for details and definitions):  

\begin{thm}[\cite{ATorelli}]
Let $S$ be any surface of infinite type. Then $\vI(S)$ is topologically generated by separating twists and bounding-pair maps. 
\end{thm} 

\subsection{Coarse boundedness}
Before we begin, we note that all the general theory about Polish groups discussed here is developed in Rosendal's forthcoming book \cite{Rosendal}. 

The theories of finitely-generated groups and locally-compact compactly-generated topological groups have many analogies, especially from the viewpoint of geometric group theory.
This is naturally due to compactness being a natural generalization of finiteness; however, as noticed by Rosendal, there is a weaker condition on topological groups that allows one to still capture many of the key aspects of the theory of locally-compact compactly-generated groups.

The key observations is to note that a compact subset of a (pseudo-)metric space always has finite diameter; it turns out this is the property to focus on.
In a Polish group \( G \), a subset \( A \) of \( G \) is \emph{coarsely bounded}\index{coarsely bounded!group}\index{group!coarsely bounded}, or \emph{CB} for short, if it has finite diameter in every continuous pseudo-metric on \( G \) (in fact, it is sufficient to only consider left-invariant continuous pseudo-metrics).
A Polish group is \emph{coarsely bounded}, or \emph{CB}, if it is coarsely bounded as a subset; it is \emph{locally coarsely bounded}\index{coarsely bounded!locally}, or \emph{locally CB}, if there exists a coarsely-bounded open neighborhood of the identity; it is \emph{CB generated}\index{coarsely bounded!generating set}\index{group!CB generated} if there exists a coarsely bounded set algebraically generating the group.

One should naturally think of CB as a generalization of compact,  locally CB as a generalization of locally compact, and CB generated as a generalization of compactly generated.
Conveniently, every CB generated Polish group is locally CB \cite[Theorem 2.30]{Rosendal} (note: it is not the case that every compactly-generated group is locally compact, e.g. \( (\QQ,+) \) is compactly generated but not locally compact).

From the point of view of this survey, the main result of the theory of CB-generated Polish groups is that, up to quasi-isometry, they have a well-defined metric.
In particular, CB-generated Polish groups have a well-defined geometry and they can be studied through the lens of geometric group theory.
Let us now describe this result.

A left-invariant continuous pseudo-metric \( d \) is \emph{maximal} if for any other left-invariant continuous pseudo-metric \( d' \) there exits constants \( K,L \geq 0 \) such that \( d' < K\cdot d +L \).
In particular, up to quasi-isometry, if a maximal pseudo-metric exists, then it is unique. 
Before stating the theorem, a subset of a Polish space is \emph{analytic} if it is the continuous image of a Polish space. 
Now, combining pieces of Theorem 1.2, Proposition 2.52, Theorem 2.53, and Example 2.54 from \cite{Rosendal}, we have:

\begin{thm}[\cite{Rosendal}]
Let \( G \) be a CB-generated Polish group. Then: 
\begin{enumerate}
\item \( G \) admits a left-invariant continuous maximal metric \( d \).
\item \( G \) has an analytic symmetric coarsely-bounded generated set; moreover, \( G \) equipped with the word metric associated to any such generating set is quasi-isometric to \( (G,d) \). 
\end{enumerate}
\end{thm}

Note that the metric topology associated to a word metric is always discrete
 and hence cannot be continuous on a non-discrete topological group. However, the
above theorem tells us that (non-continuous) word metrics capture the geometry of the group.

In recent work, Mann--Rafi \cite{MannRafi} classify the CB, locally CB, and CB-generated mapping class groups.  
The most general version of their result is a bit technical to state, so we will state a specific case that captures the main flavor.
It is a classical result of Mazurkiewicz and Sierpinski \cite{MS} that every countable compact Hausdorff topological space is homeomorphic to an ordinal space of the form \( \omega^\alpha n + 1 \), where \( \alpha \) is a countable ordinal, \( n \) is a natural number, and \( \omega \) is the first infinite ordinal.

\begin{thm}[\cite{MannRafi}]
Let \( S \) be an infinite-type surface so that either every end of \( S \) is planar or every end of \( S \) is non-planar. 
If the end space of \( S \) is countable and homeomorphic to \( \omega^\alpha n +1 \), then
\begin{enumerate}
\item \( \Map(S) \) is CB if and only if \( n = 1 \).
\item If \( n \geq 2 \) and \( \alpha \) is a successor ordinal, then \( \Map(S) \) is CB generated, but not CB.
\item If \( n \geq 2 \) and \( \alpha \) is a limit ordinal, then \( \Map(S) \) is locally CB, but not CB generated. 
\end{enumerate}
\end{thm}

The full statement of Mann--Rafi's theorem involves generalizing the trichotomy above to uncountable end spaces; they do this by introducing a partial order on the ends.
We encourage the interested reader to see their paper for details; we believe the various cases described will be essential for researchers interested in proving results about all big mapping class groups. 

For examples, the mapping class group of the Loch Ness monster surface is CB as is  the mapping class group of the flute surface.  
Also, though it does not fit into the countable version of the Mann--Rafi theorem given above, the mapping class group of the Cantor tree surface is CB. 
For \( n \in \NN \), let \( \Omega_n \) denote the infinite-genus surface with \( n \) ends, all of which are non-planar. If $n\ge 2$, then \( \Map(\Omega_n) \) is CB generated, but not CB; in particular, \( \Map(\Omega_n) \) is not quasi-isometric to \( \Map(\Omega_1) \) if $n\ge 2$. Therefore, we ask:

\begin{qu}
Are \( \Map(\Omega_n) \) and \( \Map(\Omega_m) \) quasi-isometric if and only if \( n = m \)?\index{mapping class group!quasi-isometry class}
\end{qu}

As a complementary question, we propose: 

\begin{qu}
Are there computable quasi-isometry invariants of CB-generated big mapping class groups (e.g.~geometric rank)?
\end{qu}

\subsection{Automatic continuity}

A topological group \( G \) has the \emph{automatic continuity property}\index{automatic continuity} if every  abstract group homomorphism from \( G \) to a separable topological group is continuous.  
There is a beautiful history to studying automatic continuity given in \cite{RosendalSurvey}; however, we only discuss several relevant examples (and non-examples). 

For a non-example, consider the following: the real line \( \RR \) and the real plane \( \RR^2 \), each equipped with the standard Euclidean topology and the group operation of (vector) addition, are isomorphic as groups.
To see this, observe that both \( \RR \) and \( \RR^2 \) are infinite-dimensional vector spaces over the rationals \( \QQ \) with bases of cardinality \( 2^{\aleph_0} \) and hence they are isomorphic.
However, \( \RR \) and \( \RR ^2 \) are not homeomorphic and hence this group isomorphism cannot be continuous.

For examples, none of which are trivial, the homeomorphism group\index{automatic continuity!homeomorphism group} of the Cantor set \cite{KechrisTurbulence} as well as the homeomorphism group of any closed manifold \cite{RosendalAutomatic, MannAutomatic} has the automatic continuity property. 
The automatic continuity property for homeomorphism groups (and some diffeomorphism groups) has been key to recent developments in approaches to the dimension growth question of Ghys \cite{Ghys1} regarding actions of infinite groups on compact manifolds (e.g. Chen--Mann \cite{ChenStructure}, Hurtado \cite{HurtadoContinuity}).
The application of automatic continuity in understanding the rigidity of homeomorphism groups of compact manifolds motivates us to ask about automatic continuity in mapping class groups, where there are also open rigidity questions (see Section \ref{sec:algebra}).

\begin{qu}
\label{q:acp}
Classify the surfaces \( S \) for which the groups \( \Homeo(S) \) and/or \( \Map(S) \) have the automatic continuity property.
\end{qu}

Recently, building on her previous work \cite{MannAutomatic}, Mann proved that the homeomorphism group of any manifold that can be realized as the interior of a compact manifold with boundary has the automatic continuity property \cite{MannAutomatic2}.
In the same article, Mann gave the first examples of infinite-type surfaces (e.g. the sphere minus a Cantor set) whose homeomorphism groups have the automatic continuity property.
Mann's result actually shows these groups have a stronger property (they are Steinhaus), which passes to quotients and hence yields: 

\begin{thm}[{\cite[Corollary 2.1]{MannAutomatic2}}]
Let \( S \) be an infinite-type surface of finite genus whose space of ends is of the form \( C \sqcup F \), where \( C \) is a Cantor space and \( F \) is a finite discrete space.
Then, \( \Map(S) \) has the automatic continuity property.
\index{automatic continuity!mapping class group}\index{mapping class group!automatic continuity}
\end{thm}

In \cite[Example 2.3]{MannAutomatic2}, Mann also gives an example of an infinite-type surface whose  homeomorphism group and mapping class group do not have the automatic continuity property. 

All the arguments establishing automatic continuity for the homeomorphism groups mentioned above rely on the same core technique, which unfortunately does not readily extend to non-compact surfaces with infinite-genus nor finite-genus with non-perfect end space.


\section{Algebraic aspects}
\label{sec:algebra}

\subsection{Algebraic rigidity}\index{rigidity!algebraic} \index{mapping class group!rigid}
\label{subsec:ivanov}
In this subsection, all surfaces are assumed to have empty boundary. A classical result of Ivanov \cite{Iva2} asserts that, with several well-understood exceptions, every automorphism of the mapping class group of a finite-type surface $S$ is induced by a homeomorphism of $S$.
Ivanov gave a simplified proof of this result using the curve complex in \cite{Iva}; however, in this case, he assumes the underlying surface has genus at least two.
This simplified proof was adapted to the remaining cases by Korkmaz \cite{Kor} and Luo \cite{Luo} independently.
In the infinite-type setting, the analogous result was established by Bavard--Dowdall--Rafi \cite{BDR}; namely, one has: 

\begin{thm}[\cite{BDR}]
\label{thm:BDR}
For any infinite-type surface $S$, \[\Aut(\Map(S)) \cong \Map^\pm(S).\]
\end{thm}

The idea of the proof of Theorem \ref{thm:BDR} is similar in spirit to that of Ivanov, adapted to the context of infinite-type surfaces. First, the authors prove that an element of $\Map(S)$ is supported on a finite-type subsurface of $S$ if and only if its conjugacy class is countable, and from this they obtain an algebraic characterization of Dehn twists, similar to Ivanov's original one, which is preserved by automorphisms. As a consequence, any given automorphism of $\Map(S)$ induces a simplicial automorphism of the curve complex $\vC(S)$ which in turn, by Theorem \ref{thm:ivanovcurves}, is induced by an element of $\Map^\pm(S)$. At this point, the mapping class obtained this way coincides with  the original automorphism  on every Dehn twist, from which one quickly deduces that they are equal.

\subsubsection{Injective and surjective homomorphisms} Ivanov's theorem gave rise to a large number of stronger rigidity results about mapping class groups. For instance, a result of Ivanov--McCarthy \cite{IMc}  asserts that mapping class groups of surfaces of genus at least three are {\em co-Hopfian}, that is, every injective endomorphism is an automorphism.
Hence, every injective endomorphism is induced by a homeomorphism of the underlying surface. The analog in the infinite-type setting is not known: 

\begin{qu}
Are mapping class groups of infinite-type surfaces co-Hopfian? 
\label{qu:cohopf}
\end{qu}

One of the main hurdles in this direction is that, for infinite-type surfaces, simplicial injections of the curve complex into itself need not come from mapping classes, in stark contrast to the case of finite-type surfaces (see \cite{Hernandez} for the strongest result of this type).  An example of this, for surfaces of infinite genus, may be found in \cite[Lemma 5.3]{HV}. We now present another instance of this phenomenon, which can be easily generalized to other punctured surfaces: 

\begin{ex}[Non-surjective simplicial injections between curve graphs]
Let $S$ be the flute surface.  As such, we may realize $S$ as the surface obtained by removing from $\mathbb{S}^2$ a convergent sequence together with its limit point. 

Fix a hyperbolic structure on $S$, and realize every simple closed curve on $S$ by its unique geodesic representative. Since there are only countably many simple closed curves on $S$, we may pick a point $p$ in the complement of the union of all the simple closed geodesics. Therefore we obtain a map $h:\vC(S) \to \vC(S\setminus \{p\})$ which is easily seen to be injective, since two curves that are disjoint on $S$ remain disjoint after puncturing. Finally, observe that $S\setminus \{p\}$ is homeomorphic to $S$, but that the map $h$ is not induced by a homeomorphism, as it is not surjective. This finishes the example. 
\end{ex}

With respect to surjective homomorphisms, a group is \emph{Hopfian} if every surjective endomorphism is an automorphism.  
It is an exercise to show that every finitely-generated residually-finite group is Hopfian; hence, mapping class groups of finite-type surfaces are Hopfian.
It is therefore natural to ask if big mapping class groups are Hopfian.
But, we quickly find a counterexample:

\begin{ex}[Non-Hopfian mapping class group]
Let \( E \) be a closed subset of the Cantor set such that the set \( E' \) of accumulation points of \( E \)  satisfies \( E' \neq E \) and \( E' \) is homeomorphic to \( E \).
For example, the ordinal space  \( \omega^\omega +1 \)  has this property.
Embed \( E \) into the 2-sphere \( \mathbb S^2 \). 
We then have that the embedding \( \mathbb S^2 \smallsetminus E \hookrightarrow \mathbb S^2 \smallsetminus E' \) induces a forgetful homomorphism \( \Map(\mathbb S^2 \smallsetminus E) \to \Map(\mathbb S^2 \smallsetminus E') \) that is surjective, but not injective.
Now, \( \mathbb S^2 \smallsetminus E \) is homeomorphic to \( \mathbb S^2 \smallsetminus E' \) and hence we see there exists a surjective endomorphism of \( \Map( \mathbb S^2 \smallsetminus E ) \) that fails to be an automorphism.
Note that the forgetful map exists only because \( E\smallsetminus E' \)---the set of isolated points of \( E \)---is invariant under the action of \( \Map(\mathbb S^2 \smallsetminus E) \).  
\end{ex}

\begin{qu}
If a surjective endomorphism of a mapping class group fails to be an automorphism, is it necessarily a forgetful homomorphism?
\end{qu}

\subsubsection{General homomorphisms} A result of Souto and the first author \cite{AS}  describes all non-trivial homomorphisms $\PMap(S) \to \PMap(S')$, where the genus of $S$ is at least six and the genus of $S'$ is less than twice the genus of $S$, showing that they arise as combinations of {\em subsurface inclusions}, {\em forgetting punctures}, and  {\em deleting boundary components}. 
A homomorphism between mapping class groups that comes from a {\em manipulation at the level of the underlying surfaces} is called \emph{geometric}.

Other than Theorem \ref{thm:BDR}, there are no results of this kind in the context of infinite-type surfaces. 
In fact, as a consequence of Theorem \ref{thm:APV} below, if \( S \) has at least two non-planar ends then there are non-geometric endomorphisms of \( \PMap(S) \).
However, all these examples factor through the (non-trivial) abelianization of \( \PMap(S) \).
An ambitious question is to ask if this is the only way to produce non-geometric endomorphisms:

\begin{qu}
Let $S$ be a surface of infinite type with no boundary. Does every non-geometric endomorphism of $\PMap(S)$  factor through its abelianization?  
\end{qu}

A much more humble question to which we do not know the answer (although we expect it to be negative) is: 

\begin{qu}
Let $S$ be Jacob's ladder surface and let $S'$ be the Loch Ness monster. Are there any homomorphisms $\PMap(S) \to \PMap(S')$ with non-abelian image? 
\end{qu}

\subsubsection{Rigidity of subgroups.}
In fact, the aforementioned result of Ivanov \cite{Iva} applies to injections between finite-index subgroups of mapping class groups. In other words, it asserts that the {\em abstract commensurator} $\Comm(\Map(S))$ of $\Map(S)$ is equal to $\Map^\pm(S)$, provided the genus of $S$ is large enough. For infinite-type surfaces, the analog is due to Bavard--Dowdall--Rafi \cite{BDR} (the proof is the same as for Theorem \ref{thm:BDR}): 

\begin{thm}[\cite{BDR}]
For any infinite-type surface $S$, \[\Comm(\Map(S)) \cong \Map^\pm(S).\]
\end{thm}

In \cite{ATorelli}, it is shown that $\vI(S)$ is also algebraically rigid; more concretely:

\begin{thm}[\cite{ATorelli}]
\label{thm:torelliautos}
For any infinite-type surface $S$, \[\Aut(\vI(S)) \cong \Comm(\vI(S)) \cong \Map^\pm(S).\]
\end{thm}

The equivalent statement for finite-type surfaces was proved by Farb--Ivanov \cite{FaIv} for automorphisms, and by Brendle-Margalit \cite{BeMa} for commensurations.

We remark that it is not known whether $\vI(S)$ has any finite-index subgroups at all; hence we ask: 

\begin{qu}
Does $\vI(S)$ have any proper finite-index subgroups? 
\end{qu}

Note that if the answer to the above question were negative, then $\Comm(\vI(S))$ would be equal to $\Aut(\vI(S))$ {\em a priori}.

Finally, we should mention a recent theorem of Brendle--Margalit \cite{BeMa2} (for closed surfaces) and McLeay \cite{McLeay} (for surfaces with punctures) which vastly generalizes the theorems above, proving that every normal subgroup which contains elements of {\em sufficiently small support} has the extended mapping class group as its automorphism and abstract commensurator group. In the setting of infinite-type surfaces one expects fewer necessary conditions, as the following result of McLeay \cite{McLeay2} shows:

\begin{thm}[\cite{McLeay2}]
Let $S$ be the Cantor tree surface. If $N$ is any normal subgroup of $\Map(S)$, then \[\Aut(N) \cong \Map^\pm(S).\]
\end{thm}

Though not directly a rigidity result, we finish this subsection by recalling a result of Lanier--Loving \cite{LanierLoving} that fits with the discussion:

\begin{thm}[{\cite{LanierLoving}}]
If \( S \) is an infinite-type surface, then every normal subgroup has trivial center.
\end{thm}

\subsection{Abelianization}\index{mapping class group!abelianization}

A classical result of Powell \cite{Powell}, building up on previous work of Mumford  \cite{Mum} and Birman \cite{Birman}, shows that the abelianization of the mapping class group of a closed surface of genus at least three is trivial. Moreover, the lantern relation can be used to establish the same result for all finite-type surfaces:

\begin{thm}[see {\cite[Theorem 5.2]{FM}}]
Let $S$ be a finite-type surface of genus at least $3$. Then $\PMap(S)$ has trivial abelianization. 
\end{thm}

By Proposition \ref{prop:gencompact}, $\Map_c(S)$ is a direct limit of finite-type mapping class groups, and hence: 

\begin{cor}
Let $S$ be a surface of genus at least 3. Then $\Map_c(S)$ has trivial abelianization. 
\label{cor:compactperfect}
\end{cor}

We would like to promote the above corollary to a statement about the pure mapping class group, and here is one instance where automatic continuity is incredibly useful. Indeed, a result of Dudley \cite{Dudley} asserts that if $G$ is a Polish group, then any homomorphism $G \to \ZZ$ is continuous. Combining this with Corollary \ref{cor:compactperfect}, we have:

\begin{thm}
Let $S$ be a surface of genus at least 3. Then, every homomorphism \[\overline{\Map_c(S)} \to \ZZ\] is trivial. In other words, \[H^1(\overline{\Map_c(S)},\ZZ) = \{1\}.\]
\end{thm}

In light of Theorem \ref{thm:PV} above, this has the following consequence: 

\begin{cor}
Let $S$ be a surface with at most one non-planar end. Then $H^1(\PMap(S),\ZZ) = \{1\}.$
\end{cor}

However, in \cite{APV} it was shown that the situation for general infinite-type surfaces is rather different. Namely, one has: 

\begin{thm}[\cite{APV}]
\label{thm:APV}
Let $S$ be a surface of genus at least two, and let $\hat S$ denote the result of filling every planar end of $S$. Then \[H^1(\PMap(S),\ZZ) \cong H_1^{\text{sep}}(\hat S, \ZZ),\] where the latter group is the subgroup of $H_1(\hat{S}, \ZZ)$ generated by homology classes with separating representatives. 
\index{mapping class group!cohomology}\index{group!cohomology}
\end{thm}

In particular, $H^1(\PMap(S),\ZZ)$ is not trivial as soon as $S$ has at least two non-planar ends. A natural problem is: 

\begin{prob}
Compute the low-dimensional (co-)homology groups of \( \Map(S) \) and \( \PMap(S) \).
\end{prob}

In his original blogpost, Calegari \cite{Cal} showed that the mapping class group of the Cantor tree surface is uniformly perfect, which implies that both \( H_1 \) and \( H^1 \) are trivial (with integer coefficients).
Recently, Calegari--Chen have computed the second homology; we record both results below:

\begin{thm}[{\cite{Cal, CaChen2}}]
Let \( \Gamma \) denote the mapping class group of the Cantor tree surface. Then  \( H^1(\Gamma, \ZZ) \), \( H_1(\Gamma, \ZZ) \) and \( H^2(\Gamma, \ZZ) \) are trivial, and  \( H_2(\Gamma, \ZZ) = \ZZ/2\ZZ \).
\end{thm}

The following structural result about pure mapping class groups provides the core piece in the proof of Theorem \ref{thm:APV}; compare with Theorem \ref{thm:topgen} above:

\begin{thm}
\label{cor:semidirect}
For any surface $S$, we have \[\PMap(S) = \overline{\Map_c(S)} \rtimes \prod_{s\in \vS} \langle h_s \rangle,\] where the rightmost group is a direct product of cyclic groups generated by pairwise-commuting handle shifts $h_s$, where $s$ ranges over a free basis of $H_1^{\text{sep}}(\hat{S},\ZZ)$. 
\end{thm}

Theorem \ref{thm:APV} leaves out some low-genus cases, which were subsequently settled by Domat--Plummer \cite{DP}. More concretely, they proved the following result for genus-one surfaces:

\begin{thm}[\cite{DP}]
Let $S$ be an infinite-type surface of genus one. Then \[H^1(\PMap(S),\ZZ)= 0.\]
\end{thm}

For an infinite-type surface $S$ of genus-zero  the situation is different, for in this case there is a surjective homomorphism $\PMap(S) \to \FF_2$, the free group on two generators, since the pure mapping class group of a four-times punctured sphere is isomorphic to $\FF_2$. Nevertheless, Domat--Plummer prove: 

\begin{thm}[\cite{DP}]
Let $S$ be an infinite-type surface of genus zero. Then $H^1(\PMap(S),\ZZ)$ contains uncountably many classes which do not come from forgetful maps to spheres with finitely many punctures. 
\end{thm}

\subsection{Quantifying rigidity}
In Section \ref{subsec:ivanov}, we saw that automorphisms of mapping class groups are geometric.
In fact, something stronger is true: outside several low-complexity cases, given two surfaces \( S_1 \) and \( S_2 \) any isomorphism \( \Map(S_1) \to \Map(S_2) \) (or \( \PMap(S_1) \to \PMap(S_2)) \) is induced by a homeomorphism \( S_1 \to S_2 \) (this is shown in \cite{BDR} in the infinite-type setting and can be deduced in the finite-type setting from \cite{Iva2, Kor, Luo}). 
In particular, in the finite-type setting, using the virtual cohomological dimension \cite{Harer2} and algebraic rank \cite{BLM} of \( \Map(S) \), it is possible to determine the topology of \( S \) from algebraic invariants of \( \Map(S) \).
Given that rigidity holds in big mapping class groups, it should be possible to do the same:

\begin{qu}
\label{qu:rigidity}
Is there a list of algebraic invariants of \( \Map(S) \) that determine the topology of \( S \)?
\end{qu}

Let us provide some  examples connecting algebraic invariants of \( \Map(S) \) and the  topology of \( S \).
First, we have the following corollary of Theorem \ref{thm:APV}:

\begin{cor}[{\cite{APV}}]
The algebraic rank of \( H^1(\PMap(S), \ZZ) \) is:\index{mapping class group!rigid}\index{rigidity!quantifying}
\begin{itemize}
\item 0 if and only if \( S \) has at most one non-planar end.
\item \( n \in \NN \) if and only if \( S \) has \( n+1 \) non-planar ends.
\item infinite if and only if \( S \) has infinitely many non-planar ends.
\end{itemize}
\end{cor}

Next, recall that a group is \emph{residually finite}\index{group!residually finite} if and only if the intersection of all its normal subgroups is the identity. 

\begin{thm}[{\cite{PV}}]
Let \( S \) be any surface.\index{mapping class group!residual finiteness}
\begin{itemize}
\item \( \PMap(S) \) is residually finite if and only if \( S \) has finite genus.
\item \( \Map(S) \) is residually finite if and only if \( S \) is of finite type.
\end{itemize}
\end{thm}

Now, it follows from the work of Bavard--Walker \cite{BW} that if \( S \) has an isolated planar end then \( \PMap(S) \) is circularly orderable (though not equivalent, the reader can read this as ``acts faithfully on the circle'').
Moreover, by forthcoming work of Aougab, Patel, and the second author \cite{AoPV}, every finite group can be realized as a subgroup of \( \PMap(S) \) whenever \( S \) has infinite-genus and no planar ends.
Combining these facts, with the two results mentioned in this subsection and the fact that \( \Aut(\PMap(S)) \cong \Map^\pm(S) \) when \( S \) is of infinite-type \cite{BDR}, we are able  to give a complete answer to Question \ref{qu:rigidity} for a countably infinite family of surfaces:

\begin{thm}[{\cite{AoPV}}]
For \( n \in \NN \), let \( \Omega_n \) denote the \( n \)-ended infinite-genus surface with no planar ends
and let \( G = \PMap(S) \) for some surface \( S \).
The surface \( S \) is homeomorphic to \( \Omega_n \) if and only if  \( G \) satisfies each of the following properties:
\begin{enumerate}
\item \( G \) is not residually finite,
\item \( G \) is not circularly orderable,
\item \( H^1(G, \ZZ) \) has rank \( n -1 \), and
\item \( G \) is finite index in \( \Aut(G) \).
\end{enumerate}
\end{thm}

\subsection{Homology representation}\index{mapping class group!homology representation}
As mentioned in \ref{subsec:ivanov}, there is a homomorphism \[\rho_S: \Map(S) \to \Aut(H_1(S,\ZZ)),\] given by the action of mapping classes on the homology of the surface. For finite-type surfaces with at most one puncture or boundary component, the {\em algebraic intersection} pairing of homology classes is a symplectic form, and one shows that the homomorphism \[\Map(S) \to \Sp(2g,\ZZ),\] where $g$ is the genus of $S$, is surjective; see \cite[Section 6]{FM} for details. 

The homology representation for infinite-type surfaces has been studied by Fanoni, Hensel, and the second author \cite{FHV}.
In the infinite-type setting, there is only one surface with at most one end, namely the Loch Ness monster surface; in this case, it turns out an analogous result holds: 

\begin{thm}[\cite{FHV}]
Let $S$ be the Loch Ness monster surface. Then the image of the homology representation is the subgroup of $\Aut(H_1(S,\ZZ))$ consisting of those elements which preserve the algebraic intersection form. In other words, \[\Im(\phi_S) = \Sp(\NN,\ZZ).\]
\end{thm}

For surfaces with more than one end (or boundary component), preserving algebraic intersection is not enough to characterize the image of \( \rho_S \) in \( H_1(S, \ZZ) \) (this is true in both the finite-type and infinite-type settings).
In the same article \cite{FHV}, the authors give a characterization of the image of \( \rho_S \) for an arbitrary surface \( S \) in terms of preserving a filtration of the first homology.
The full statement is a bit technical, so we refer the interested reader directly to \cite{FHV}.

\subsection{Nielsen realization}\index{Nielsen realization} Kerckhoff's {\em Nielsen Realization Theorem} \cite{Ker} asserts that every finite subgroup of the mapping class group of a finite-type surface $S$ of negative Euler characteristic lifts to $\Homeo(S)$; moreover, it may be realized as a subgroup of the isometry group of some hyperbolic metric on $S$. 

In the context of big mapping class groups, the analogous statement has been obtained by Afton--Calegari--Chen--Lyman \cite{AfLy}: 

\begin{thm}
Let $S$ be a surface of infinite type. Then every finite subgroup of $\Map(S)$ lifts to $\Homeo^+(S)$.
Moreover, every finite group can be realized as a group of isometries of some hyperbolic metric on \( S \).
\index{Nielsen realization!infinite-type surfaces}
\end{thm} 

We should also note that there is analog of Nielsen realization in the setting of analytically-infinite Riemann surfaces due to Markovic \cite{Markovic}.
A \emph{hyperbolic Riemann surface} is a complex 1-manifold whose universal cover is isomorphic to the unit disk. 

\begin{thm}[{\cite{Markovic}}]
Let \( S \) be an infinite-type surface and let \( G \) be a subgroup of \( \Map(S) \).
If there exists a hyperbolic Riemann surface \( X \) homeomorphic to \( S \) and a constant \( K > 1 \) such that every element of \( G \) can be realized by a \( K \)-quasi-conformal homeomorphism \( X \to X \), then there is a hyperbolic Riemann surface \( Y \) such that \( Y \) is quasi-conformally equivalent to \( X \) and \( G < \rm{Isom}(Y) \).
\end{thm}

%

\subsection{The relation with Thompson groups}\index{group!Thompson's}\index{Thompson groups} Thompson's groups $F$, $T$ and $V$ constitute prominent examples of discrete subgroups of $\Homeo(C)$, the homeomorphism group of the Cantor set. Among many other features, they are infinite groups of type $F_\infty$, and which have simple commutator subgroup (in fact, $V$ itself is simple). We now briefly review the construction of these groups, referring the reader to the standard reference \cite{CFP} for a thorough treatment of Thompson's groups. 

\subsubsection{Thompson's groups} Let $\mathcal T$ be a rooted binary tree, noting that its space of ends of $\mathcal T$ is homeomorphic to the Cantor set $C$. The tree $\mathcal T$ has a natural left-to-right orientation once we fix a realization of $\mathcal{T}$ as a subset of the hyperbolic plane. With respect to this orientation, given a subtree of $\mathcal T$ with $n$ leaves, we may order its set of leaves using the numbers $1, \ldots, n$, so that the numbers increase from left to right.

 Let $\tau, \tau'$ be subtrees of $\mathcal T$ with the same number of leaves, and such that both contain the root of $\mathcal T$. If $\sigma$ is a bijection between the sets of leaves of $\tau$ and $\tau'$, then the triple $(\tau,\tau',\sigma)$ extends in a natural way to a homeomorphism of $C$. Of course, the same homeomorphism may be induced by different such triples (obtained by {\em expanding} and {\em contracting} a given finite subtree), and Thompson's group $V$ is the group of equivalence classes of such triples. In turn, Thompson's group $T$ (resp.~$F$) corresponds to the case when the bijection $\sigma$ is a cycle (resp.~the identity).

\subsubsection{Asymptotic mapping class groups}\index{mapping class group!asymptotic} We now explain the relation between Thompson's groups and big mapping class groups. To this end, let $S$ denote either the Cantor tree surface or the blooming Cantor tree surface.
 In these particular cases, the exact sequence \eqref{eq:SEShomeo} reads \begin{equation}
1 \to \PMap(S) \to \Map(S) \to \Homeo(C) \to 1.
\label{eq:SESCantor}
\end{equation}
Over the last two decades, numerous authors have given geometric constructions of finitely-presented subgroups $H$ of $\Map(S)$ for which the sequence \eqref{eq:SESCantor} restricts to 
\begin{equation}
1 \to \Map_c(S) \to H \to G \to 1,
\label{SESgenB}
\end{equation}
where $G$ is one of Thompson's groups $F$, $T$ or $V$ (or their commutator subgroups).

To the best of our knowledge, the first step in this direction was the paper of Greenberg--Sergiescu \cite{GreSe}, whose objective  was to construct an acyclic extension of $F'$, the commutator subgroup of $F$, by the braid group $B_\infty$ on infinitely many strands. This was later generalized simultaneously by Brin \cite{Brin} and Dehornoy \cite{Dehornoy} to the construction of an extension of $V$ by $B_\infty$, the so-called {\em braided Thompson groups}\index{Thompson groups!braided}. Funar--Kapoudjian \cite{FK,FK2}, and later Funar and the first author \cite{AF}, constructed finitely-generated (and often finitely-presented) extensions of $V$ by a direct limit of mapping class groups of compact surfaces. Part of the motivation \cite{FK2} is to construct a finitely-presented group whose homology agrees with 
 the {\em stable homology} of pure mapping class groups, after a seminal result of Harer \cite{Harer}.

A common feature of all of the above constructions is that they may be expressed in terms of groups of homeomorphisms of an infinite-type surface which {\em eventually} preserve some topological data; these are the {\em asymptotic mapping class groups} introduced by Funar--Kapoudjian in \cite{FK}. We now briefly recall their definition in the simpler case of a surface of genus zero.

\subsubsection{The case of the Cantor tree surface} Let $S$ be the Cantor tree surface, that is, a sphere with a Cantor set removed. Fix, once and for all, a pants decomposition $P$ of $S$ and a set $A$ of pairwise-disjoint, properly-embedded arcs on $S$ such that $S\smallsetminus A$ has exactly two connected components $\nu^\pm$, and each connected component of $S\smallsetminus P$ is intersected by exactly three elements of $A$; see  Figure \ref{fig:genus0} . The triple $(P,A, \nu^+)$ is called a {\em rigid structure}\index{surface!rigid structure on} on $S$. 
 \\
 \begin{figure}[htb]
\begin{center}
\includegraphics[width=3in,height=2in]{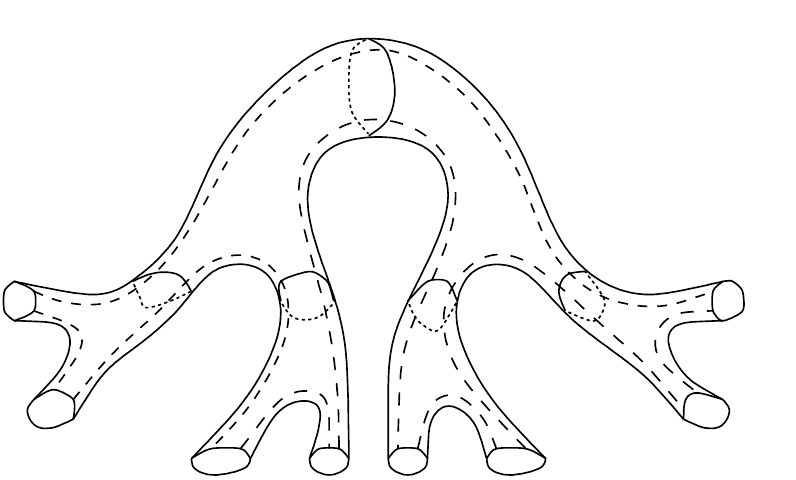} \caption{The rigid structure on $S$.} \label{fig:genus0}
\end{center}
\end{figure}
We say that a homeomorphism $f: S \to S$ is {\em asymptotically rigid}\index{homeomorphism!asymptotically rigid} if there exists a compact subsurface $X \subset S$ with $\partial X \subset P$, such that $\partial f(X) \subset P$ and the restriction homeomorphism \[f: S \smallsetminus X \to S \smallsetminus f(X)\] setwise preserves (the relevant part of) the rigid structure. The group $\mathcal{B}$ is then defined as the subgroup of $\Map(S)$ whose elements have an asymptotically rigid homeomorphism. In their paper \cite{FK}, Funar and Kapoudjian showed that the restriction of the sequence \eqref{SESgenB} yields 
\begin{equation}
1 \to \Map_c(S) \to \mathcal{B} \to V \to 1,
\label{eq:SESB}
\end{equation}
As such, $\mathcal{B}$ contains the mapping class group of every compact surface of genus zero with non-empty boundary. In light of this, the main result of \cite{FK} is rather striking: 

\begin{thm}[\cite{FK}]
The group $\mathcal{B}$ is finitely presented. 
\end{thm}

Moreover, they observed: 

\begin{prop}
The short exact sequence \eqref{eq:SESB} splits over Thompson's group $T$. As a consequence, $\mathcal{B}$ is not linear and does not have Kazhdan's Property (T).  
\end{prop}

We remark that a well-known question about finite-type mapping class groups asks whether they are linear or have Kazhdan's Property (T). 

\subsubsection{Other compact surfaces with a Cantor set removed}
The construction of asymptotic mapping class groups makes sense for arbitrary surfaces. In fact, as commented in \cite{FK}, the group constructed by Brin \cite{Brin} and Dehornoy \cite{Dehornoy} are asymptotic mapping class groups of a closed disc with a Cantor set removed, and as such embeds as a subgroup of $\mathcal{B}$. In addition, Funar and the first author \cite{AF} generalized the construction of $\mathcal{B}$ to the surface $\Sigma_g$ obtained by removing a Cantor set from a closed surface $S_g$ of finite genus $g\ge 1$. Roughly speaking, a {\em rigid structure} on $\Sigma_g$ is determined by a simple closed curve $\alpha\subset \Sigma_g$ that cuts off a once-punctured surface of genus $g$, together with a rigid structure for the planar component of $\Sigma_g$. One then defines the notion of an asymptotically rigid homeomorphism in an analogous way, and constructs the {\em asymptotic mapping class group} $\mathcal{B}_g$\index{mapping class group!asymptotic} as the subgroup of $\Map(\Sigma_g)$ whose elements have an asymptotically rigid representative. In this case, the restriction of the short exact sequence \eqref{eq:SESCantor} to the group $\vB_g$ reads
\begin{equation}
1 \to \Map_c(\Sigma_g) \to \vB_g \to V \to 1;
\end{equation}
in particular, $\vB_g$ contains the mapping class group of every compact surface of genus at most $g$ and with non-empty boundary. The following is one of the main results of \cite{AF}: 

\begin{thm}[\cite{AF}] For every $g \ge 1$, the group $\vB_g$ is finitely presented. In addition, it is not linear and does not have Kazhdan's Property (T). 
\label{thm:Bfg}
\end{thm}

In addition, in \cite{AF} the authors explore the structure of the groups $\vB_g$ in connection with mapping class groups of finite-type surfaces. For instance,  every automorphism of $\vB_g$ is induced by a homeomorphism of $\Sigma_g$ (compare with Theorems \ref{thm:BDR} and \ref{thm:torelliautos}). 

\subsubsection{The case of the blooming Cantor tree} In \cite{FK2}, Funar and Kapoudjian constructed an asymptotic mapping class group $\vB_\infty$ for the blooming Cantor tree, which we denote by $\Sigma_\infty$. In a similar fashion, the short exact sequence \eqref{eq:SESCantor}, when restricted to $\vB_\infty$, yields: 
\begin{equation}
1 \to \Map_c(\Sigma_\infty) \to \vB_\infty \to V \to 1. 
\end{equation}

The following is the main result of \cite{FK2}:

\begin{thm}[\cite{FK2}]
The group $\vB_\infty$ is finitely generated. Moreover, its rational cohomology coincides with the stable rational cohomology of the mapping class group. 
\end{thm}

Note that, while asymptotic mapping class groups of finite genus are finitely presented, the group $\vB_\infty$ is only known to be finitely generated. In light of this, we ask: 

\begin{qu}
Determine whether the asymptotic mapping class groups $\vB_n$, for $n\in \mathbb{N} \cup \{\infty\}$, satisfy stronger finiteness properties. Are they ${\rm F}_\infty$? 
\end{qu}

A positive answer to the above question, in the case of $n=0$, is conjectured in \cite[p. 967]{FK}. The question of whether $\vB_\infty$ is finitely presented appears in \cite{FKS}.

\subsubsection{A dense asymptotic mapping class group}
In addition, in \cite{AF} the authors considered a subgroup $\mathcal{H}_g$ with $\vB_g<\vH_g <\Map(\Sigma_g)$. The definition of $\vH_g$ is similar to that of $\vB_g$, without the requirement that its elements preserve the connected component $\nu^+$ appearing in the definition of rigid structure. In short, the difference between $\vB_g$ and $\vH_g$ is that the latter contains {\em half-twists} about separating curves cutting off a disk minus a Cantor set. For this reason, the group $\vH_g$ is referred to as the {\em group of half-twists}. 

A large part of the motivation for considering $\vH_g$ comes from the study of smooth mapping class groups, as explained in \cite{FNe}. Indeed, put a differentiable structure on the closed surface $S_g$ of genus $g$, and realize $C$ as the the middle-third Cantor set on a smoothly-embedded interval on $S_g$. Let $\Mod^s(S_g, C)$ denote the {\em smooth mapping class group} of the pair $(S_g,C)$, namely the group of isotopy classes of smooth diffeomorphisms of $S_g$ preserving globally the Cantor set $C$. The following is a recent result of Funar and Neretin \cite{FNe}:

\begin{thm}[\cite{FNe}, Cor. 2]
For every $g \ge 0$, we have $\vH_g \cong \Mod^s(S_g,C)$.
\end{thm}

Using the same techniques as with $\mathcal{B}_g$, Funar and the first author \cite{AF} proved:

\begin{thm}[\cite{AF}] For every $g \ge 1$, the group $\vH_g$ is finitely presented. In addition, it is not linear and does not have Kazhdan's Property (T). 
\label{thm:Bgenfg}
\end{thm}

However, a nice extra feature of the group $\vH_g$ is the following result, which should be compared with Theorem \ref{thm:PV}: 

\begin{thm}
For every $g\ge 0$ the group $\vH_g$ is dense in $\Map(\Sigma_g)$.\index{mapping class group!topologically generating}
\end{thm}

Finally, the restriction to $\vH_g$ of the sequence \eqref{eq:SESB} reads
\begin{equation}
1 \to \Map_c(\Sigma_g) \to \vH_g \to V_2[\mathbb{Z}_2] \to 1,
\end{equation}
where $V_2[\mathbb{Z}_2]$ is the {\em Higman--Thompson group}\index{group!Higman--Thompson} $V_2[\mathbb{Z}_2]$ \cite{BDJ}. A surprising result of Bleak--Donoven--Jonu\v{s}as \cite{BDJ} establishes that $V$ and $V_2[\mathbb{Z}_2]$ are conjugate as subgroups of $\Homeo(C)$ through an explicit homeomorphism of $C$ (a {\em cellular automaton}). An obvious questions then is: 

\begin{qu}
Are the groups $\vB_g$ and $\vH_g$ isomorphic? 
\end{qu}

If would be surprising if the question above had a positive answer, since isomorphisms between (sufficiently rich) subgroups of mapping class groups tend to come from surface homeomorphisms. 

We end this section with the following vague question: 

\begin{qu}
Are there other  geometrically-defined subgroups of $\Map(\Sigma_g)$ which surject to other interesting classes of subgroups of $\Homeo(C)$, such as Higman-Thompson groups, Neretin groups, etc? 
\end{qu}


\section{Geometric aspects} 
\label{sec:geometry}

Mapping class groups of finite-type surfaces have been successfully studied through their action on various combinatorial  complexes, notably the curve graph; a first instance of this is Ivanov's Rigidity Theorem mentioned in Section \ref{subsec:ivanov}. Moreover, it turns out that the geometric structure of $\vC(S)$, equipped with its natural path metric, sheds intense light on the algebraic and geometric structure of $\Map(S)$. In this direction, the following is a seminal theorem of Masur--Minsky \cite{MM1}: 

\begin{thm}[\cite{MM1}]
\label{thm:MM}
Let $S$ be a finite-type surface. If $\vC(S)$ is connected, then it is hyperbolic (in the sense of Gromov).
\end{thm}

A number of authors have proved analogous results for other combinatorial complexes associated to surfaces, such as the disk graph \cite{MS}, the non-separating curve graph ${\rm NonSep}(S)$ \cite{MS,Ham},  the arc graph $\vA(S)$ \cite{HPW}, etc. In fact, a surprising phenomenon is that the hyperbolicity constant in Theorem \ref{thm:MM}, as well as the those of other complexes, turn out to be independent of the underlying surface; we say that the corresponding family of complexes are {\em uniformly hyperbolic}\index{curve graph!uniform hyperbolicity}. The following theorem is a combination of the results of \cite{HPW,A,Bow,CRS,Ras}: 
\begin{thm}  Let \( S \) be a finite-type surface.
\begin{enumerate}
\item (\cite{HPW}) $\vA(S)$ is uniformly hyperbolic.
\item (\cite{HPW,A,Bow,CRS}) $\vC(S)$ is uniformly hyperbolic
\item (\cite{Ras}) For fixed $g$, the graph  $\rm{NonSep}(S_{g,n}))$ is hyperbolic with respect to a constant which does not depend on $n$. 
\end{enumerate}
\label{thm:unifhyp}
\end{thm}

The above result may be regarded as a curiosity at first, but it happens to be of central importance in the study of big mapping class groups, as we will explain next.

\subsection{Complexes for infinite-type surfaces} As in the finite-type case, one may be tempted to use interesting geometric properties of analogous combinatorial models, built from arcs and/or curves, in order to study mapping class groups. This initial surge of enthusiasm is thwarted by the following immediate observation; before we state it, we recall that, for an infinite-type surface $S$, the arc graph $\vA(S)$ is defined to be the simplicial graph whose vertices are properly embedded arcs on $S$ which join two (not necessarily distinct) planar ends of $S$, and where adjacency corresponds to disjointness. 

\begin{fact}Let $S$ be a surface of infinite type. Then $\vC(S)$ has diameter two. Furthermore, if $S$ has infinitely many planar ends, then $\vA(S)$ also has diameter two.
\end{fact}

However, as mentioned in the introduction, in \cite{Cal} Calegari proposed studying $\Map(\RR^2 \setminus C)$ via its action on a certain subgraph of $\mathcal{A}(\RR^2 \setminus C)$; observe that, by the above, $\vA(\RR^2 \setminus C)$ itself has diameter two. 
Calegari's idea was to consider the subgraph \( \vA_\infty \) of \( \vA(\RR^2 \smallsetminus C) \) consisting of arcs with at most one endpoint in \( C \) (hence, necessarily one end of an arc in \( \vA_\infty \) is contained in the unique isolated planar end of \( \RR^2 \smallsetminus C \)). 
The next result was proved by Juliette Bavard \cite{Bav,Bav2} proving a conjecture posed by Calegari: 

\begin{thm}
\label{thm:bavard}
$\vA_\infty$ is a Gromov-hyperbolic space of infinite diameter. 
\end{thm}

Based on this result, and with a lot of extra work, she also proved that
$\Map(\RR^2 \setminus C)$ has an infinite-dimensional space of quasi-morphisms. 
This is in stark contrast to \( \Map(\mathbb S^2 \smallsetminus C) \), which Calegari shows admits no quasi-morphisms (and even stronger, we know \( \Map(\mathbb S^2 \smallsetminus C) \) is CB \cite{MannRafi}). 
We note that the automorphism group of \( \vA_\infty \) and related graphs are computed in \cite{Sch} and shown to be the extended mapping class group.

The above theorem may be regarded as part of a more general phenomenon, which we now explain. In order to do so, we need the following terminology due to Schleimer \cite{Saulnotes}. Given a graph $\mathfrak{X}(S)$ built from arcs and/or curves on $S$, say that a subsurface $Y \subset S$ is a {\em witness} for $\mathfrak{X}(S)$ if every vertex of $\mathfrak{X}$ intersects $Y$ non-trivially. For instance, the only non-trivial witness for $\vC(S)$ is $S$ itself while, in the case of $\vA(S)$, any subsurface $Y \subset S$ which contains every puncture of $S$ is a witness.

The following theorem is a reformulated version of \cite[Theorem 1]{AV} (see also \cite[Section 6]{DFV} for another formulation). In an intuitive way, it encapsulates the idea of {\em taking a limit of a family of uniformly hyperbolic spaces}:

\begin{thm}
Let $\mathfrak{X}(S)$ be a connected $\Map(S)$-invariant graph, whose vertices are defined by finite sets of arcs or curves on $S$, and where edges correspond to bounded intersection number. Given a subsurface $Y \subset S$, define $\mathfrak{X}(Y)$ to be the full subgraph of $\mathfrak{X}(S)$ spanned by those vertices which are entirely contained in $Y$ and equip $\mathfrak X(Y)$ with the induced path metric. Suppose that: 

\begin{enumerate}
\item For every triangle $T$ in $\mathfrak{X}(S)$ there exists a finite-type witness $Y$ such that $T$ is contained in $\mathfrak{X}(Y)$ and \( \mathfrak{X}(Y) \) is connected;
\item There exists constants \( \delta, K , C > 0 \) such that for every finite-type witness $Y$ of $S$ with \( \mathfrak{X}(Y) \)  connected, the following conditions are satisfied: 
\begin{enumerate}

\item  $\mathfrak{X}(Y)$ is a \( \delta \)-hyperbolic graph of infinite diameter. 
\item The inclusion map $\mathfrak{X}(Y) \hookrightarrow \mathfrak{X}(S)$ is a \( (K,C) \)-quasi-isometric embedding.
\end{enumerate}
\end{enumerate} 
Then $\mathfrak{X}(S)$ is hyperbolic  and has infinite diameter. 

\end{thm}

Given a finite set $P$ of isolated planar ends of $S$, denote by $\vA(S;P)$ the subgraph of $\mathcal A(S)$ spanned by those arcs which have at least one endpoint in $P$; observe that every subsurface of $S$ which contains $P$ is a witness for  $\vA(S;P)$. 
The above result and the uniform hyperbolicity presented in Theorem \ref{thm:unifhyp} are used to prove the following:

\begin{thm} Let $S$ be an infinite-type surface. 
\begin{enumerate}
\item (\cite{Bav,AFP,AV}) Let $P$ be a non-empty finite set of isolated punctures of $S$. Then, $\vA(S;P)$ is hyperbolic. 
\item (\cite{Ras})If $S$ has finite genus at least $2$, then the graph ${\rm{NonSep}}(S)$ is hyperbolic. 
\end{enumerate}
\label{thm:AFP}
\index{curve graph!analogs for infinite-type surfaces}
\end{thm}

\begin{rmk}
There is a subtlety  about Theorem \ref{thm:AFP} which is worth mentioning at this point; see also \cite[Theorem 1]{AV}. Let $P, Q$ be  two  finite sets of isolated punctures of $S$, with $P \cap Q = \emptyset$, and consider the subgraph $\mathcal{A}(S; P, Q)$ of $\mathcal A(S)$ which have one endpoint in $P$ {\em and} one endpoint in $Q$. Then $\mathcal{A}(S; P, Q)$ is not hyperbolic. 

Indeed, this is a manifestation of Schleimer's {\em Disjoint Witness Property} \cite{Saulnotes}, which asserts that if a graph or curves/arcs has two disjoint witnesses of infinite diameter then it is not hyperbolic, for one may use subsurface projections to construct a quasi-isometrically embedded copy of $\mathbb{Z}^2$ inside the graph. 

Finally, observe that the graph $\mathcal{A}(S; P, Q)$ contains two disjoint witnesses, since one can take two finite-type surfaces, one containing $P$ and the other containing $Q$. This finishes the remark.
\end{rmk}

These different phenomena were clarified in subsequent work of Durham, Fanoni and the second author \cite{DFV}.
The motivation of their work was to find actions of big mapping class group not relying on isolated planar ends.
Before explaining their result, we need some definitions. 

Let $\vQ$ be a collection of pairwise-disjoint closed subsets of $\Ends(S)$.
Every separating curve on \( S \) partitions \( \Ends(S) \); let \( \rm{Sep}_2(S,\vQ) \) denote the subgraph of \( \mathcal C(S) \) consisting of separating curves on \( S \) that partition \( \vQ \) into two sets, each of cardinality at least 2 (there is a slight modification if \( |\vQ|=4 \), see \cite{DFV} for details). 

\begin{thm}[\cite{DFV}] Let $S$ be an infinite-type surface.
Let $\vQ$ be a collection of pairwise-disjoint closed subsets of $\Ends(S)$ such that, for every $\omega \in \vQ$ and every $f\in \Map(S)$, there exists $\omega' \in \vQ$  with $f(\omega) = \omega'$. 
Then, \( \rm{Sep}_2(S, \vQ) \) is hyperbolic, infinite diameter, \( \Map(S) \)-invariant, and 
 there are infinitely many mapping classes which act with positive translation length on \( \rm{Sep}_2(S, \vQ) \). 
\label{thm:DFV}
\end{thm}

For example, if \( S = \Omega_n \) (the $n$-ended infinite-genus surface with no planar ends) with \( n \geq 4 \), then \( \vQ = \Ends(S) \) satisfies the hypothesis of the above theorem. 

We note that in the days this survey was being finalized, Fanoni--Ghaswala--\mc Leay \cite{FGM} constructed new examples of hyperbolic infinite-diameter graphs that admit actions of big mapping class groups with unbounded orbits.
We direct the reader to their article for details. 

Klarreich \cite{Kla} showed that the Gromov boundary of the curve graph is $\Map(S)$-equivariantly homeomorphic to the space of {\em ending laminations} on the surface; see also \cite{Hambound} for a different argument, and Pho-On's thesis \cite{Pho-On} for an effective proof of this using the {\em unicorn} machinery of \cite{HPW}. In unpublished work, Schleimer proved that the boundary of the arc graph is naturally identified with the space of all ending laminations supported on witnesses of $S$; this is also carried out in an effective manner in Pho-On's thesis \cite{Pho-On}.

In light of these results, we ask the following natural question: 

\begin{qu}
Describe the Gromov boundary of the various hyperbolic complexes associated to an infinite-type surface $S$, ideally in terms of laminations/foliations on $S$. 
\end{qu}

For the case when the surface is $\mathbb{R}^2 \setminus C$, the Gromov boundary of the relative arc graph $\mathcal A_\infty$ of Theorem \ref{thm:bavard} is described by Bavard-Walker \cite{BW} in terms of rays on the surface.
Rasmussen \cite{Ras3} has recently reproved a result of Hamenst\"adt computing the Gromov boundary of the graph of non-separating curves and points out that his techniques can be extended to the infinite-type setting; however, the issue is a lack of understanding of laminations on infinite-type surfaces.  
We should note at this point that \v{S}ari\'c \cite{Saric} recently developed the theory of train tracks for infinite-type surfaces, which should aid in investigating laminations. 

The natural motivation for understanding the Gromov boundary is to gain insight into a potential classification of big mapping classes akin to that of the Nielsen--Thurston classification.  We should note that there is much research in this direction for quasi-conformal mapping class groups and their action on Teichm\"uller space. 

\subsection{Weak proper discontinuity and acylindricity} Let $G$ be a group acting by isometries on a hyperbolic metric space $(X,d)$. We say that the action is {\em acylindrical} if, for every $D \ge 0$, there exists $R\ge 0$ such that, for every $x,y \in X$ with $d(x,y) \ge D$, the cardinality of the set  \[\{g\in G \mid d(x,gx),d(y,gy) \le R\}\] is finite. To exclude uninteresting pathologies, we restrict our attention to actions where there are infinitely many points on the Gromov boundary of $X$ that are accumulation points of $G$-orbits; call such an action {\em non-elementary}. We say that a group is {\em acylindrically hyperbolic}\index{group!acylindrically hyperbolic} if it admits a non-elementary acylindrical action on some Gromov-hyperbolic space.

A result of Bowditch \cite{Bow-acy} asserts that, if $S$ has finite type, the action of $\Map(S)$ on $\vC(S)$ is acylindrical. Bavard--Genevois \cite{BaGe} proved that the analogous statement  does not hold for infinite-type surfaces: 

\begin{thm}[\cite{BaGe}]
If $S$ has infinite type, then $\Map(S)$ is not acylindrically hyperbolic. 
\end{thm}

Prior to the notion of acylindricity, Bestvina--Fujiwara \cite{BeF} introduced the concept of {\em weak proper discontinutity}\index{weak proper discontinuity} (WPD, for short), and used it to show that if a group has an interesting WPD action then it has an infinite-dimensional space of quasimorphisms; equivalently, its second bounded cohomology group is infinite-dimensional. We briefly recall these ideas. Let $G$ be a group acting on a Gromov-hyperbolic metric space $(X,d)$, and $g\in G$ be a loxodromic element. We say that $g$ is a {\em WPD element} if, for every $x\in X$ and every $R \ge 0$, there is $N \in \NN$ such that the set \[\{h \in G \mid d(x,h(x)) \le R, d(g^N(x),hg^N(x) \le R\}\] is finite. 
Bestvina--Fujiwara \cite{BeF} proved that, for a finite-type surface $S$, any pseudo-Anosov element of $\Map(S)$ is WPD with respect to the natural action on the curve complex. This notion was further weakened by Bromberg--Bestvina--Fujiwara \cite{BBF} to that of a {\em WWPD action}: suppose again $G$ acts on a hyperbolic space $X$, and let $g$ be a loxodromic element of $G$ with fixed points $\eta^\pm$ on the Gromov boundary $\partial X$ of $X$. We say that $g$ is a {\em WWPD} element if, for every sequence $\{h_n\}_{n\in \NN}$ of elements of $G$, with $h_n(\eta^+) \to \eta^+$ and $h_n(\eta)^- \to \eta^-$, there exists $N \in \NN$ such that, for all $n \ge N$, one has \[ h_n(g^+)=g^+ \hspace{.5cm} \text{ and } \hspace{.5cm}  h_n(g^-)=g^-.\]

The existence of WWPD elements\index{mapping class group!WWPD elements} of big mapping class groups has been recently studied by Rasmussen \cite{Ra2}. 
Let $S$ be an infinite-type surface with at least one isolated puncture $p$, and let $\vA(S,p)$ be the relative arc graph of $S$ based at $p$. Rasmussen proved: 

\begin{thm}[{\cite{Ra2}}]
An element $g\in \Map(S)$ is WWPD with respect to the action of $\Map(S)$ on $\vA(S,p)$ if and only if there exists a finite-type $g$-invariant subsurface $Y \subset S$, with $p \in Y$,  such that  the restriction of $g$ to $Y$ is pseudo-Anosov. 
\end{thm}

As a consequence, he deduces that a  class of subgroups of $\Map(S)$ have infinite-dimensional second bounded cohomology. 

We finish with mentioning a very recent construction of Morales--Valdez \cite{MV}, in which they construct examples of mapping classes which act loxodromically on $\vA_\infty$ and do not preserve any finite-type subsurface.


\printindex

\end{document}